\documentclass[12pt]{article}
\usepackage{amssymb,amsmath,stmaryrd}
\usepackage{hyperref}
\hypersetup{pdfpagemode=FullScreen,colorlinks=true}
\def\R{\mathbb{R}}
\def\Z{\mathbb{Z}}
\def\N{\mathbb{N}}

\def\eps{\epsilon}
\newtheorem{thm}{Theorem}
\newtheorem{prop}{Proposition}
\newtheorem{cor}[prop]{Corollary}
\newtheorem{lem}[prop]{Lemma}

\newtheorem{rem}[prop]{Remark}
\newtheorem{ex}[prop]{Example}
\newtheorem{defi}[prop]{Definition}
\def\n#1{|\hskip-1.5pt|#1|\hskip-1.5pt|}
\newenvironment{pf}{\begin{trivlist}\item[]{\bf Proof\ }}
{\mbox{}\hfill\rule{.08in}{.08in}\end{trivlist}}

\title{Large scale conformal maps}

\author{Pierre Pansu\footnote{P.P. is supported by MAnET Marie Curie
Initial Training Network, by Agence Nationale de la Recherche, ANR-10-BLAN 116-01 GGAA and ANR-15-CE40-0018 SRGI. P.P. gratefully acknowledges the hospitality of Isaac Newton Institute, of EPSRC under grant EP/K032208/1, and of Simons Foundation.
}}

\begin{document}

\maketitle

\begin{center}
Titre fran\c cais : Applications conformes \`a grande \'echelle
\end{center}

\abstract{Roughly speaking, let us say that a map between metric spaces is \emph{large scale conformal} if it maps packings by large balls to large quasi-balls with limited overlaps. This quasi-isometry invariant notion makes sense for finitely generated groups. Inspired by work by Benjamini and Schramm, we show that under such maps, some kind of dimension increases: exponent of polynomial volume growth for nilpotent groups, conformal dimension of the ideal boundary for hyperbolic groups. A purely metric space notion of $\ell^p$-cohomology plays a key role.}

\medskip
\begin{center}
\textbf{R\'esum\'e}
\end{center}
	Grosso modo, une application entre espaces m\'etriques est \emph{conforme \`a grande \'echelle} si elle envoie tout empilement de grandes boules sur une collection de grandes quasi-boules qui ne se chevauchent pas trop. Cette notion est un invariant de quasi-isom\'etrie, elle s'\'etend aux groupes de type fini. En s'inspirant de travaux de Benjamini et Schramm, on montre qu'en pr\'esence d'une telle application, une sorte de dimension doit augmenter : il s'agit de l'exposant de croissance polyn\^omiale du volume pour les groupes nilpotents, de la dimension conforme du bord pour les groupes hyperboliques. Une nouvelle d\'efinition, purement m\'etrique, de la cohomologie $\ell^p$ joue un r\^ole important.

\tableofcontents

\section{Introduction}

\subsection{Microscopic conformality}

Examples of conformal mappings arose pretty early in history: the stereographic projection, which is used in astrolabes, was known to ancient Greece. The metric distorsion of a conformal mapping can be pretty large. For instance, the Mercator planisphere (1569) is a conformal mapping of the surface of a sphere with opposite poles removed onto an infinite cylinder. Its metric distorsion (Lipschitz constant) blows up near the poles, as everybody knows. Nevertheless, in many circumstances, it is possible to estimate metric distorsion, and this lead in the last century to metric space analogues of conformal mappings, known as quasi-symmetric or quasi-M\"obius maps. Grosso modo, these are homeomorphisms which map balls to quasi-balls. Quasi means that the image $f(B)$ is jammed between two concentric balls, $B\subset f(B)\subset \ell B$, with a uniform $\ell$, independent of location or radius of $B$.

\subsection{Mesoscopic conformality}

Some evidence that conformality may manifest itself in a discontinuous space shows up with Koebe's 1931 circle packing theorem. A circle packing of the 2-sphere is a collection of interior-disjoint disks. The incidence graph of the packing has one vertex for each circle and an edge between vertices whenever corresponding circles touch. Koebe's theorem states that every triangulation of the 2-sphere is the incidence graph of a disk packing, unique up to M\"obius transformations. Thurston conjectured that triangulating a planar domain $\Omega$ with a portion of the incidence graph of the standard equilateral disk packing, and applying Koebe's theorem to it, one would get a numerical approximation to Riemann's conformal mapping of $\Omega$ to the round disk. This was proven by Rodin and Sullivan, \cite{RS}, in 1987. This leads us to interpret Koebe's circle packing theorem as a mesoscopic analogue of Riemann's conformal mapping theorem.

\subsection{A new class of maps}

In this paper, we propose to go one step further and define a class of large scale conformal maps. Roughly speaking, large scale means that our definitions are unaffected by local changes in metric or topology. Technically, it means that pre- or post-composition of large scale conformal maps with quasi-isometries are again large scale conformal. This allows to transfer some techniques and results of conformal geometry to discrete spaces like finitely generated groups, for instance.

\subsection{Examples}

In first approximation, a map between metric spaces is large scale conformal if it maps every packing by sufficiently large balls to a collection of large quasi-balls which can be split into the union of boundedly many packings. We postpone till next section the rather technical formal definition. Here are a few sources of examples.
\begin{itemize}
  \item Quasi-isometric embeddings are large scale conformal.
  \item Snowflaking (i.e. replacing a metric by a power of it) is large scale conformal.
  \item Power maps $z\mapsto z|z|^{K-1}$ are large scale conformal for $K\geq 1$. They are not quasi-isometric, nor even coarse embeddings.
  \item Compositions of large scale conformal maps are large scale conformal.
\end{itemize}
For instance, every nilpotent Lie or finitely generated group can be large scale conformally embedded in Euclidean space of sufficiently high dimension, \cite{A}. Every hyperbolic group can be large scale conformally embedded in hyperbolic space of sufficiently high dimension, \cite{BoS}.

\subsection{Results}

Our first main result is that a kind of dimension increases under large scale conformal maps. The relevant notion depends on classes of groups.

\begin{thm}
\label{dd}
If $G$ is a finitely generated or Lie nilpotent group, set $d_1(G)=d_2(G)=$ the exponent of volume growth of $G$. If $G$ is a finitely generated or Lie hyperbolic group, let $d_1(G):=CohDim(G)$ be the infimal $p$ such that the $\ell^p$-cohomology of $G$ does not vanish. Let $d_2(G):=ConfDim(\partial G)$ be the Ahlfors-regular conformal dimension of the ideal boundary of $G$.

Let $G$ and $G'$ be nilpotent or hyperbolic groups. If there exists a large scale conformal map $G\to G'$, then $d_1(G)\leq d_2(G')$.
\end{thm}

Theorem \ref{dd} is a large scale version of a result of Benjamini and Schramm, \cite{BS}, concerning packings in $\R^d$. The proof follows the same general lines but differs in details. The result is not quite sharp in the hyperbolic group case, since it may happen that $d_1(G) < d_2(G)$, \cite{BP}. However equality $d_1(G) = d_2(G)$ holds for Lie groups, their lattices and also for a few other finitely generated examples.

\medskip

Our second result is akin to the fact that maps between geodesic metric spaces which are uniform/coarse embeddings in both directions must be quasi-isometries.

\begin{thm}
\label{isomorphism}
Let $X$ and $X'$ be bounded geometry manifolds or polyhedra. Assume that $X$ and $X'$ have isoperimetric dimension $>1$. Every homeomorphism $f:X\to X'$ such that $f$ and $f^{-1}$ are large scale conformal is a quasi-isometry.
\end{thm}

Isoperimetric dimension is defined in subsection \ref{isop}. Examples of manifolds or polyhedra with isoperimetric dimension $>1$ include universal coverings of compact manifolds or finite polyhedra whose fundamental group is not virtually cyclic.

This is a large scale version of the fact that every quasiconformal diffeomorphism of hyperbolic space is a quasi-isometry. This classical result generalizes to Riemannian $n$-manifolds whose isoperimetric dimension is $>n$. The new feature of the large scale version is that isoperimetric dimension $>1$ suffices.

\subsection{Proof of Theorem \ref{dd}}

Our main tool is a metric space avatar of energy of functions, and, more generally, norms on cocycles giving rise to $L^p$-cohomology. Whe\-reas, on Riemannian $n$-manifolds, only $n$-energy $\int|\nabla u|^n$ is conformally invariant, all $p$-energies turn out to be large scale conformal invariants. Again, we postpone the rather technical definitions to Section \ref{energy} and merely give a rough sketch of the arguments.

Say a locally compact metric space is $p$-parabolic if for all (or some) point $o$, there exist compactly supported functions taking value 1 at $o$, of arbitrarily small $p$-energy. For instance, a nilpotent Lie or finitely generated group $G$ is $p$-parabolic iff $p \geq d_1(G)$. Non-elementary hyperbolic groups are never $p$-parabolic. If $X$ has a large scale conformal embedding into $X'$ and $X'$ is $p$-parabolic, so is $X$. This proves Theorem 1 for nilpotent targets.

If a metric space has vanishing $L^p$-cohomology, maps with finite $p$-energy have a limit at infinity. We show that a hyperbolic group $G'$ admits plenty of functions with finite $p$-energy when $p>ConfDim(G')$. First, such functions separate points of the ideal boundary $\partial G'$ (this result also appears in \cite{BnK}). Second, for each ideal boundary point $\xi$, there exists a finite $p$-energy function with a pole at $\xi$. The first fact implies that, if $CohDim(G)>ConfDim(G')$, any large scale conformal map $G\to G'$ converges to some ideal boundary point $\xi$, the second leads to a contradiction. This proves Theorem 1 for hyperbolic domains and targets.

In \cite{BS}, Benjamini and Schramm observed that vanishing of reduced cohomology suffices for the previous argument to work, provided the domain is not $p$-parabolic. This proves Theorem 1 for nilpotent domains and hyperbolic targets.

In the hyperbolic to hyperbolic case, one expects $CohDim$ to be replaced by $ConfDim$. For this, one could try to reconstruct the ideal boundary of a hyperbolic group merely in terms of finite $p$-energy functions, in the spirit of \cite{R} and \cite{B}.

\subsection{Proof of Theorem \ref{isomorphism}}

Following H. Gr\"otzsch, \cite{grotzsch}, we define $1$-capacities of compact sets $K$ in a locally compact metric space $X$ by minimizing $1$-energies of compactly supported functions taking value 1 on $K$. Then we minimize capacities of compact connected sets joining a given pair of points to get a pseudo-distance $\delta$ on $X$. This is invariant under homeomorphisms which are large scale conformal in both directions. If $X$ has bounded geometry, $\delta$ is finite. If $X$ has isoperimetric dimension $d>1$, then $\delta$ tends to infinity with distance. This shows that homeomorphisms which are large scale conformal maps in both directions are coarse embeddings in both directions, hence quasi-isometries. It turns out that all finitely generated groups have isoperimetric dimension $>1$, but virtually cyclic ones. 

\subsection{Larger classes}
\label{larger}

Some of our results extend to wider classes of maps. If we merely require that balls of a given range of sizes are mapped to quasi-balls which are not too small, we get the class of \emph{uniformly conformal maps}. It is stable under precomposition with arbitrary uniform (also known as coarse) embeddings. We can show that no such map can exist between nilpotent or hyperbolic groups unless the inequalities of Theorem \ref{dd} hold. 

\begin{cor}\label{udd}
We keep the notation of Theorem \ref{dd}. Let $G$ and $G'$ be hyperbolic groups. If there exists a uniform embedding $G\to G'$, then $d_1(G)\leq d_2(G')$.
\end{cor}
We note that special instances of this Corollary have been obtained by D. Hume, J. Mackay and R. Tessera by a different method, \cite{HMT}. Their results apply in particular to M. Bourdon's rich class of (isometry groups of) Fuchsian buildings, see section \ref{lack}.

If we give up the restriction on the size of the images of balls, we get the even wider class of \emph{coarse conformal maps}. It is stable under post-composition with quasi-symmetric homeomorphisms. New examples arise, such as stereographic projections, or the Poincar\'e model of hyperbolic space and its generalizations to arbitrary hyperbolic groups. However, when targets are smooth, coarse conformal maps are automatically uniformly conformal, hence similar results hold.

\begin{cor}\label{cdd}
No coarse conformal map can exist between a finitely generated group $G$ and a nilpotent Lie group $G'$ equipped with a Riemannian metric unless $G$ is itself virtually nilpotent, and $d_1(G)\leq d_1(G')$. Also, no coarse conformal map can exist from a hyperbolic group $G$ to a bounded geometry manifold quasiisometric to a hyperbolic group $G'$ unless $d_1(G)\leq d_2(G')$.
\end{cor}

\subsection{Organization of the paper}

Section \ref{notions} contains definitions, basic properties and examples of coarse, uniform, rough and large scale conformal maps. In Section \ref{qsstructure}, the notion of a quasi-symmetry structure is introduced, as a tool to handle hyperbolic metric spaces: every such space has a rough conformal map onto a product of quasi-metric quasi-symmetry spaces, as shown in Section \ref{poinc}. The existence of this map has the effect of translating large scale problems into microscopic analytic issues. The definition of energy in Section \ref{energy} comes with moduli of curve families and parabolicity. It culminates with the proof that several families of quasi-symmetry spaces are parabolic. $L^p$-cohomology of metric spaces is defined in Section \ref{lp}, where the main results relating parabolicity, $L^p$-cohomology and coarse conformal maps are proven. Section \ref{lack} draws consequences for nilpotent or hyperbolic groups, concluding the proof of Theorem \ref{dd}. The material for the proof of Theorem \ref{isomorphism} is collected in Section \ref{iso}. As a byproduct, we find conditions on a pair of spaces $X,X'$ in order that coarse conformal maps $X\to X'$ be automatically uniformly conformal, this provides the generalizations selected in Corollary \ref{cdd}.

\subsection{Acknowledgements}

The present work originates from a discussion with James Lee and Itai Benjamini on conformal changes of metrics on graphs, cf. \cite{ST}, during the Institut Henri Poincar\'e trimester on Metric geometry, algorithms and groups, see \cite{metric2011}. It owes a lot to Benjamini and Schramm's paper \cite{BS}. The focus on the category of metric spaces and large scale conformal maps was triggered by a remark by Jonas Kahn.

\section{Coarse notions of conformality}
\label{notions}

A sphere packing in a metric space $Y$ is a collection of interior-disjoint balls. The incidence graph $X$ of the packing has one vertex for each ball and an edge between vertices whenever corresponding balls touch. A packing may be considered as a map from the vertex set $X$ to $Y$, that maps the tautological packing of $X$ (by balls of radius $1/2$) to the studied packing of $Y$.

We modify the notion of a sphere packing in order to make it more flexible. In the domain, we allow radii of balls to vary in some finite interval $[R,S]$. In the range, we replace collections of disjoint balls with collections of balls with bounded multiplicity (unions of boundedly many packings). We furthermore insist that $\ell$-times larger concentric balls still form a bounded multiplicity packing.

The resulting notion is invariant under coarse embeddings between domains and quasi-symmetric maps between ranges. Therefore, it is a one-sided large scale concept (in terms of domain, not of range). It is reminiscent of conformality since it requires that spheres be (roughly) mapped to spheres. Whence the term ``coarsely conformal''. In order to get a class which is invariant under post-composition with quasi-isometries, we shall introduce a subclass of ``large scale conformal'' maps.

\subsection{Coarse, uniform, rough and large scale conformality}

Let $X$ be a metric space. A ball $B$ in $X$ is the data of a point $x\in X$ and a radius $r\geq 0$. For brevity, we also denote $B(x,r)$ by $B$. If $\lambda\geq 0$, $\lambda B$ denotes $B(x,\lambda r)$. For $S\geq R\geq 0$, let $\mathcal{B}^{X}_{R,S}$ denote the set of balls whose radius $r$ satisfies $R\leq r\leq S$.

\begin{defi}
Let $X$ be a metric space. An \emph{$(\ell,R,S)$-packing} is a collection of balls $\{B_j\}$, each with radius between $R$ and $S$, such that the concentric balls $\ell B_j$ are pairwise disjoint. An $(N,\ell,R,S)$-packing is the union of at most $N$ $(\ell,R,S)$-packings. 
\end{defi}
The balls of an $(N,\ell,R,S)$-packing, $N\geq 2$, are not disjoint (I apologize for this distorted use of the word \emph{packing}), but no more than $N$ can contain a given point.

\begin{defi}
Let $X$ and $X'$ be metric spaces. Let $f:X\to X'$ be a map. Say $f$ is \emph{$(R,S,R',S')$-coarsely conformal} if there exists a map  
$$B\mapsto B',\quad\mathcal{B}^{X}_{R,S}\to \mathcal{B}^{X'}_{R',S'},$$
and, for all $\ell'\geq 1$, an $\ell\geq 1$ and an $N'$ such that 
\begin{enumerate}
  \item For all $B\in\mathcal{B}^{X}_{R,S}$, $f(B)\subset B'$.
  \item If $\{B_j\}$ is a $(\ell,R,S)$-packing of $X$, then $\{B'_j\}$ is an $(N',\ell',R',S')$-packing of $X'$.
\end{enumerate}
\end{defi}

\begin{defi}
Let $f:X\to X'$ be a map between metric spaces.
\begin{enumerate}
  \item We say that $f$ is \emph{coarsely conformal} if there exists $R>0$ such that for all finite $S\geq R$, $f$ is $(R,S,0,\infty)$-coarsely conformal.
    \item We say that $f$ is \emph{uniformly conformal} if for every $R'>0$, there exists $R>0$ such that for all finite $S\geq R$, $f$ is $(R,S,R',\infty)$-coarsely conformal.
  \item We say that $f$ is \emph{roughly conformal} if there exists $R>0$ such that $f$ is $(R,\infty,0,\infty)$-coarsely conformal.
  \item We say that $f$ is \emph{large scale conformal} if for every $R'>0$, there exists $R>0$ such that $f$ is $(R,\infty,R',\infty)$-coarsely conformal.
\end{enumerate}
\end{defi}

Here is the motivation for these many notions. The main technical step in our theorems applies to the larger class of coarse conformal maps. To turn this into a large scale notion, one needs to forbid the occurrence of small balls, whence the slightly more restrictive uniform variant, to which all our results apply. Uniformly conformal maps do not form a category, it is the smaller class of large scale conformal maps which does. The Poincar\'e models of hyperbolic metric spaces are crucial tools, but these maps are not large scale conformal, since small balls do occur in the range, merely roughly conformal.

\begin{prop}\label{composition}
The four classes enjoy the following properties
\begin{itemize}
  \item Large scale conformal $\implies$ roughly conformal $\implies$ coarsely conformal.
  \item Large scale conformal $\implies$ uniformly conformal $\implies$ coarsely conformal.
  \item Let $X$, $X'$ and $X''$ be metric spaces. Let $f:X\to X'$ be $(R,S,R',S')$-coarsely conformal. Let $f':X'\to X''$ be $(R',S',R'',$ $S'')$-coarsely conformal. Then $f'\circ f : X\to X''$ is $(R,S,R'',S'')$-coarsely conformal. 
  \item Large scale conformal maps can be composed. Large scale conformal maps up to \emph{translations}, i.e. self-maps that move points a bounded distance away, constitute the morphisms of the \emph{large scale conformal category}.
  \item Roughly conformal maps can be precomposed with large scale conformal maps. Precomposing a roughly conformal map with a uniformly conformal map yields a coarsely conformal map.
  \item Uniformly conformal maps between locally compact metric spaces are automatically proper. 
\end{itemize}
\end{prop}

\begin{pf} 
Composing maps
\begin{align*}
\mathcal{B}^{X}_{R,S}\to \mathcal{B}^{X'}_{R',S'}\to \mathcal{B}^{X''}_{R'',S''},
\end{align*}
we get, for every ball $B$ in $X$, balls $B'$ in $X'$ and $B''$ in $X''$ such that $f(B)\subset B'$, $f'(B')\subset B''$, hence $f'\circ f(B)\subset B''$. Furthermore, we get, for every $\ell''\geq 1$, a scaling factor $\ell'$ and a multiplicity $N''$, and then a scaling factor $\ell$ and a multiplicity $N'$. Given an $(\ell,R,S)$-packing of $X$, the corresponding balls can be split into at most $N'$ $(\ell',R',S')$-packings of $X'$. For each sub-packing, the corresponding balls in $X''$ can be split into at most $N''$ $(\ell'',R'',S'')$-packings of $X''$. This yields a total of at most $N'N''$ $(\ell'',R'',S'')$-packings of $X''$, i.e. a $(N'N'',\ell'',R'',S'')$-packing, as desired.

Properness of uniformly conformal maps is proven by contradiction. If $f:X\to X'$ is uniformly conformal but not proper, there exists a sequence $x_j\in X$ such that $f(x_j)$ has a limit $x'\in X'$. Fix $R'>0$ and $\ell'\geq 1$, get $R>0$, $\ell\geq 1$ and $N'$. One may assume that $d(x_j,x_{j'})>2\ell R$ for all $j'\not=j$. Then $\{B(x_j,R)\}$ is a $(\ell,R,R)$-packing. There exist balls $B'_j\supset f(B(x_j,R))$ which form a $(N',\ell',R',\infty)$-packing. For $j$ large, $x'\in B'_j$, contradicting multiplicity $\leq N'$. One concludes that $f$ is proper.
\end{pf}

In the next subsections, we shall relate our large scale conformal definitions with classical notions and collect examples.

\subsection{Quasi-symmetric maps}

Notions of quasi-conformal maps on metric spaces have a long history, see \cite{HK}, \cite{He}, \cite{Ty}.

\begin{ex}
By definition, a homeomorphism $f:X\to X'$ is \emph{quasi-symmetric} if there exists a homeomorphism $\eta:\R_+ \to\R_+$ such that for every triple $x,y,z$ of distinct points of $X$,
  \begin{align*}
\frac{d(f(x),f(y))}{d(f(x),f(z))}\leq \eta(\frac{d(x,y)}{d(x,z)}).
\end{align*}
\end{ex}

\begin{prop}\label{qs=>cc}
Quasi-symmetric homeomorphisms are roughly (and thus coarsely) conformal.
\end{prop}

\begin{pf}
When $B=B(x,r)$, we define $B'$ to be the smallest ball centered at $f(x)$ which contains $f(B)$. Let $\rho'$ be its radius. Let $y\in B$ be such that $d(f(x),f(y))=\rho'$. If $z\in f^{-1}(\ell' B')$, $d(f(x),f(z))\leq \ell'\rho'$, so
\begin{align*}
\frac{d(f(x),f(y))}{d(f(x),f(z))}\geq \frac{\rho'}{\ell'\rho'}=\frac{1}{\ell'}.
\end{align*}
By quasi-symmetry, this implies that $\eta(\frac{d(x,y)}{d(x,z)})\geq \frac{1}{\ell'}$, and thus $d(x,z)\leq \frac{1}{\eta^{-1}(\frac{1}{\ell'})}d(x,y)$. In other words, $z\in \ell B$ with
$\ell=\frac{1}{\eta^{-1}(\frac{1}{\ell'})}$. We conclude that, for every $\ell'\geq 1$, there exist $\ell\geq 1$ and a correspondance $B\mapsto B'$ such that $f(B)\subset B'$ and $f^{-1}(\ell'B')\subset\ell B$. 

If $B_1$ and $B_2$ are balls in $X$ such that $\ell B_1$ and $\ell B_2$ are disjoint, then $f^{-1}(\ell' B'_1)\subset \ell B_1$ and $f^{-1}(\ell' B'_2)\subset \ell B_2$ are disjoint as well, hence $\ell'B'_1$ and $\ell'B'_2$ are disjoint, thus $f$ is $(0,\infty,0,\infty)$-coarsely conformal. A fortiori, $f$ is $(R,\infty,0,\infty)$-coarsely conformal for $R>0$, so $f$ is roughly conformal.
\end{pf}

\medskip

Note that $R$ and $S$ play no role when checking that quasi-symmetric maps are coarsely conformal, and no guarantee on radii of balls in the range is given (i.e. such maps are $(0,\infty,0,\infty)$-coarsely conformal). In a sense, for the wider class of maps we are interested in, three of the quasi-symmetry assumptions are relaxed:
\begin{itemize}
  \item Maps need not be homeomorphisms.
  \item The quasi-symmetry estimate applies only to balls in a certain range $[R,S]$ of radii.
  \item Centers need not be mapped to centers.
\end{itemize}

\begin{prop}\label{qc=>cc}
Globally defined quasiconformal mappings of Euclidean space $\R^n$, $n\geq 2$, are quasi-symmetric, hence roughly and coarsely conformal. 
\end{prop}

\begin{pf}
Although not explicitely stated, this follows from the proof of Theorem 22.3, page 78, in \cite{vaisala}. In $\R^n$, there is a uniform lower bound $h(\frac{b}{a})>0$ for the conformal capacity of condensers $(C_0,C_1)$ such that $C_0$ connects 0 to the $a$-sphere and $C_1$ connects the $b$-sphere to infinity. Let $f:\R^n\to\R^n$ be $K$-quasiconformal and map 0 to 0. Let $l=\min\{|f(x)|\,;\,|x|=a\}$ and $L=\max\{|f(x)|\,;\,|x|=b\}$, let $C_0=f^{-1}(B(0,l))$ and $C_1=f^{-1}(\R^n\setminus B(0,L))$. Then $cap_n(C_0,C_1)\geq h(\frac{b}{a})$. On the other hand, 
\begin{align*}
cap_n(C_0,C_1)\leq K\,cap_n(f(C_0,C_1))=K\omega_{n-1}\,\log(\frac{L}{l})^{1-n},
\end{align*}
this yields an upper bound on $\frac{L}{l}$ in terms of $\frac{b}{a}$, proving that $f$ is quasi-symmetric.
\end{pf}

\begin{ex}
\label{power}
For all $K>0$, the map $z\mapsto z|z|^{K-1}$ on $\R^n$ is roughly conformal. If $K\geq 1$, it sends large balls to large balls, hence it is large scale conformal.
\end{ex}
If $n\geq 2$, Proposition \ref{qc=>cc} applies. Its 1-dimensional analogue $f:x\mapsto x|x|^{K-1}$ is roughly conformal as well, and large scale conformal if $K\geq 1$. This can be checked directly. 

\subsection{Coarse quasi-symmetry}

From the proof of Proposition \ref{qs=>cc}, we see that the following is a sufficient condition for coarse conformality.

\begin{defi}
Let $X$ and $X'$ be metric spaces. Let $f:X\to X'$ be a map. Say $f$ is $(R,S,R',S')$-\emph{coarsely quasi-symmetric} if it comes with a map 
$$B\mapsto B',\quad\mathcal{B}^{X}_{R,S}\to \mathcal{B}^{X'}_{R',S'},$$
and for all $\ell'\geq 1$, there exist $\ell\geq 1$ such that
\begin{enumerate}
  \item for all $B\in\mathcal{B}^{X}_{R,S}$, $f(B)\subset B'$.
  \item If $\ell B_1\cap\ell B_2=\emptyset$, then $\ell'B'_1\cap\ell' B'_2=\emptyset$.
\end{enumerate}

Say that $f$ is \emph{coarsely quasi-symmetric} if there exists $R>0$ such that for all $S\geq R $, $f$ is $(R,S,0,\infty)$-coarsely quasi-symmetric.
\end{defi}

\begin{lem}
Let $X$ and $X'$ be metric spaces. Let $f:X\to X'$ be a $(R,S,R',S')$-coarsely quasi-symmetric map. Then $f$ is $(R,S,R',S')$-coarsely conformal. In particular, coarsely quasi-symmetric maps are coarsely conformal.
\end{lem}

\begin{pf}
Coarse quasi-symmetry amounts to requiring $N'=1$ in the definition of coarse conformality.
\end{pf}

\begin{rem}
If a map $f:X\to X'$ is $(R,S,R',S')$-coarsely quasi-symmetric, then the ball correspondence $B\mapsto B'$ satisfies, for all $\ell'\geq 1$, $f^{-1}(\ell' B')\subset \ell B$ for a suitable $\ell$.
\end{rem}

\begin{pf}
If $x\notin 2\ell B$, the two element collection $\{B_1=B,B_2=B(x,R)\}$ is a $(\ell,R,S)$-packing. Thus $\{B',B'_2\}$ is a $(\ell',R',S')$-packing. Therefore $f(x)\notin \ell'B'$, $x\notin f^{-1}(\ell' B')$. This shows that $f^{-1}(\ell' B')\subset 2\ell B$.
\end{pf}

\subsection{Quasi-M\"obius maps}

Let $X$ be a metric space. If $a,b,c,d\in X$ are distinct, their \emph{cross-ratio} is
\begin{align*}
[a,b,c,d]=\frac{d(a,c)}{d(a,d)}\frac{d(b,d)}{d(b,c)}.
\end{align*}
An embedding $f:X\to Y$ is \emph{quasi-M\"obius} if there exists a homeomorphism $\eta:\R_+ \to \R_+$ such that for all quadruples of distinct points $a,b,c,d\in X$,
\begin{align*}
[f(a),f(b),f(c),f(d)]\leq\eta([a,b,c,d]).
\end{align*}
Note that if $f$ is a homeomorphism, $f^{-1}$ is quasi-M\"obius as well. 

The main examples are 
\begin{itemize}
  \item inversions $x\mapsto \frac{x}{|x|^2}$ in Banach spaces,
  \item the stereographic projection $\R^n \to S^n$, 
  \item its complex, quaternionic and octonionic versions, sometimes known as \emph{Cayley transforms}, where $\R^n$ is replaced with a Heisenberg group equipped with a Carnot-Carath\'eodory metric.
  \item the action of a hyperbolic group on its ideal boundary is (uniformly) quasi-M\"obius.
\end{itemize}

According to J. V\"ais\"al\"a, \cite{V}, if $X$ and $Y$ are bounded, quasi-M\"obius maps $X\to Y$ are quasi-symmetric. If $X$ and $Y$ are unbounded, a quasi-M\"obius map $X\to Y$ that tends to infinity at infinity is quasi-symmetric. We note an other situation where a quasi-M\"obius map is coarsely conformal.

\begin{lem}
\label{qm=>cc}
Assume $X$ is unbounded and $Y$ is bounded. Let $f:X\to Y$ be a quasi-M\"obius embedding. Assume that $f$ has a limit at infinity. Then $f$ is coarsely conformal.
\end{lem}

\begin{pf}
Fix $S>0$. Let $y=\lim_{x\to\infty} f(x)$. Fix some origin $o\in X$ et let $o'=f(o)$. Consider a triple $a,b,c\in X$ such that $d(c,o)\geq 4S$, $d(c,a)\leq S$ and $d(c,b)\leq S$. Then $d(a,o)\geq 3S$, $d(b,o)\geq 3S$ and $d(a,b)\leq 2S$, thus 
\begin{align*}
\frac{1}{3}\leq 1-\frac{2S}{d(b,o)}\leq\frac{d(a,o)}{d(b,o)}\leq 1+\frac{2S}{d(b,o)}\leq 2,
\end{align*}
Let $x\in X$ satisfy $d(c,x)\geq 4S$. Then $d(x,a)\geq 3S$ and $d(x,b)\geq 3S$ as well, so
\begin{align*}
\frac{1}{6}\leq [a,b,x,o]\leq 6.
\end{align*}
It follows that $\frac{1}{\eta^{-1}(\frac{1}{6})}\leq [f(a),f(b),f(x),f(o)]\leq \eta(6)$. As $x$ tends to $\infty$, this cross-ratio tends to
\begin{align*}
[f(a),f(b),y,o']=\frac{d(f(a),y)}{d(f(a),o')}\frac{d(f(b),o')}{d(f(b),y)},
\end{align*}
showing that $\frac{d(f(a),y)}{d(f(b),y)}$ stays bounded above and below.

Since
\begin{align*}
[a,b,c,\infty]=\frac{d(a,c)}{d(b,c)},\quad [f(a),f(b),f(c),y]=\frac{d(f(a),f(c))}{d(f(b),f(c))}\frac{d(f(b),y)}{d(f(a),y)},
\end{align*}
the ratio $\frac{d(f(a),f(c))}{d(f(b),f(c))}$ is bounded above and below in terms of $\frac{d(a,c)}{d(b,c)}$ only. The argument in the previous proof show that $f$ is $(0,S,0,\infty)$-large scale quasi-symmetric. One concludes that $f$ is coarsely conformal. \end{pf}

\subsection{The Cayley transform}

The stereographic projection (or Cayley transform) extends to all metric spaces. It is specially well suited to the class of Ahlfors regular metric spaces.

Recall that a metric space $X$ is \emph{$Q$-Ahlfors regular at scale $S$} if it admits a measure $\mu$ and a constant $C(S)$ such that
\begin{align*}
\frac{1}{C}\,r^Q\leq\mu(B(x,r))\leq C\,r^Q
\end{align*}
for all $x\in X$ and $r\leq S$. We abbreviate it in $Q$-Ahlfors regular if $S=$ diameter$(X)$. Locally $Q$-Ahlfors regular means $Q$-Ahlfors regular at all scales (with constants depending on scale).

\begin{ex}
$\R$, intervals of $\R$ and $\R/\Z$ are 1-Ahlfors regular. Carnot groups are Ahlfors regular. Snowflaking, i.e. replacing the distance $d$ by $d^\alpha$ for some $0<\alpha<1$, turns a $Q$-Ahlfors regular space into a $\frac{Q}{\alpha}$-Ahlfors regular space. The product of $Q$- and $Q'$-Ahlfors regular spaces is a $Q+Q'$-Ahlfors regular space. The ideal boundary of a hyperbolic group equipped with a visual quasi-metric is Ahlfors regular, \cite{S}, \cite{Coo}.
\end{ex}

\begin{lem}[V\"ais\"al\"a, \cite{V}]
\label{arqm}
Every metric space $X$ has a quasi-M\"o\-bius embedding into a bounded metric space $\dot{X}$. If $X$ is $Q$-Ahlfors regular, so is $\dot{X}$.
\end{lem}

\begin{pf}
Use the Kuratowski embedding $X\to L$, where $L=L^{\infty}(X)$. Then embed $L$ into $L\times\R$ and apply an inversion. This is a quasi-M\"obius map onto its image, which is bounded and homeomorphic to the one-point completion $\dot{X}$ of $X$. 

If $X$ is $Q$-Ahlfors regular with Hausdorff measure $\mu$, let $\nu$ be the measure with density $|\dot{x}|^{2Q}$ with respect to the pushed forward measure on $\dot{X}$. Since the inverse image of $B(\dot{x},r)$ is roughly equal to $B(x,r|\dot{x}|^{-2})$, where $|\dot{x}|=|x|^{-1}$,  
\begin{align*}
\nu(B(\dot{x},r))\sim |\dot{x}|^{2Q}\mu(B(x,r|\dot{x}|^{-2}))\sim\mu(B(x,r))\sim r^{Q}.
\end{align*}
This shows that $\dot{X}$ is $Q$-Ahlfors regular as well. \end{pf}

\subsection{Uniform/coarse embeddings}

\begin{defi}
A map $f:X\to X'$ is a \emph{uniform} or \emph{coarse embedding} if for every $T>0$ there exists $\tilde{T}$ such that for every $x_1$, $x_2\in X$,
\begin{align*}
d(x_1,x_2)\leq T &\implies d(f(x_1),f(x_2))\leq \tilde{T},\\ d(f(x_1),f(x_2))\leq T &\implies d(x_1,x_2)\leq \tilde{T}.
\end{align*}\end{defi}

\begin{lem}
\label{coarse=>lsc}
For every $0<R\leq S<+\infty$, a uniforme (or coarse) embedding is $(R,S,R',S')$-large scale quasi-symmetric for some positive and finite $R'$ and $S'$, hence $(R,S,R',S')$-coarsely conformal. In particular, uniform (or coarse) embeddings are uniformly conformal.
\end{lem}

\begin{pf}
Assume that $f:X\to X'$ is a coarse embedding controlled by function $T\mapsto \tilde{T}$. Pick an arbitrary $R>0$ and an arbitrary $S\geq R$. Set $R'=S'=\tilde{S}$, $U=\ell' \tilde{S}$ and $\ell=\frac{\tilde{U}}{R}$.

When $B=B(x,r)$, $r\in[R,S]$, we define $B'=B(f(x),\tilde{S})$. Then $f(B)\subset B'$.

On the other hand, if $x'\in f^{-1}(\ell'B')$, then $d(f(x'),f(x))\leq\ell'\tilde{S}=U$, thus $d(x',x)\leq\tilde{U}$. This shows that $f^{-1}(\ell'B')\subset \ell B$ with $\ell=\frac{\tilde{U}}{R}$. In other words, $f$ is $(R,S,\tilde{S},\tilde{S})$-large scale quasi-symmetric.
\end{pf}

\begin{rem}
Conversely, if the metric spaces are geodesic, every $(R,S,$ $R',S')$-coarse conformal map with $0<R\leq S<+\infty$ and $0<R'\leq S'<+\infty$ is a coarse embedding.

We see that $(R,S,R',S')$-coarse conformality does not bring anything new while $R,R'>0$ and $S,S'<\infty$, at least in the geodesic world. Similarly, classical conformal mappings with Jacobian bounded above and below are bi-Lipschitz. Conformal geometry begins when $S'=\infty$ or $R'=0$.
\end{rem}
Indeed, let $f:X\to X'$ be a $(R,S,R',S')$-coarse conformal map. Then $R$-balls are mapped into $S'$-balls. If $X$ is geodesic, this implies that 
$$d(f(x_1),f(x_2))\leq \frac{S'}{R}d(x_1,x_2)+S'.$$ 
Conversely, set $\ell'=(T+2S')/2R'$ and let $\ell$ be the corresponding scaling factor in the domain. Let $B_1=B(x_1,S)$ and $B_2=B(x_2,S)$. If $d(x_1,x_2)>2\ell S$, then $\ell B_1$ and $\ell B_2$ are disjoint. Let $B'_1=B(x'_1,r_1)\supset f(B_1)$ and $B'_2=B(x'_2,r_2)\supset f(B_1)$ be the corresponding balls in $X'$. Since $\{B_1,B_2\}$ is a $(1,\ell,R,S)$-packing, $\{B'_1,B'_2\}$ is a $(1,\ell',R',S'\}$ packing, so both $R'\leq r_i\leq S'$ and $\ell B'_1\cap \ell B'_2=\emptyset$, hence $d(x'_1,x'_2)> \ell'r_1+\ell'r_2\geq 2\ell'R'$. Since both $d(f(x_i),x'_i)\leq S'$, $d(f(x_1),f(x_2))>2\ell'R'-2S'\geq T$.

\begin{lem}
\label{lincoarse=>lsc}
Quasi-isometric embeddings are large scale conformal. 
\end{lem}

\begin{pf}
The assumption means that
\begin{align*}
\frac{1}{L}\,d(x_1,x_2)-C\leq d(f(x_1),f(x_2))\leq L\,d(x_1,x_2)+C.
\end{align*}
To a ball $B=B(x,r)$ in $X$, we attach $B'=B(f(x),Lr+C)$. Given balls $B_1=B(x_1,r_1)$ and $B_2=B(x_2,r_2)$, assume that concentric balls $\ell B_1$ and $\ell B_2$ are disjoint. Then, looking at an (almost) minimizing path joining $x_1$ to $x_2$, one sees that $d(x_1,x_2)\geq \ell(r_1+r_2)$. This implies
\begin{align*}
d(f(x_1),f(x_2))\geq \frac{1}{L}\ell(r_1+r_2)-C \geq \ell'(r_1+r_2)
\end{align*}
provided $r_1+r_2 \geq 2R$ and $\ell'\leq\frac{\ell}{L}-\frac{C}{2R}$. If so, the concentric balls $\ell'B'_1$ and $\ell'B'_2$ are disjoint. Therefore, for every $\ell'\geq 1$, there exists $\ell\geq 1$ such that $f$ maps $(1,\ell,R,\infty)$-packings to $(1,\ell',\frac{R}{L}-C,\infty)$-packings, and thus $(\ell,R,\infty)$-packings to $(\ell',\frac{R}{L}-C,\infty)$-packings for all $N$.\end{pf}

\begin{cor}
Quasi-isometries between geodesic metric spaces are large scale conformal. 
\end{cor}
Indeed, quasi-isometries between geodesic metric spaces are controlled by affine functions in both directions.

\begin{ex}
Orbital maps of injective homomorphisms between isometry groups of locally compact metric spaces are coarse embeddings. They are rarely quasi-isometric.
\end{ex}
For instance, the control function $T\mapsto \tilde{T}$ of the horospherical embedding of $\R^n$ into $H^{n+1}$ is logarithmic. Given a Euclidean ball $B$ of radius $R$, the hyperbolic ball $B'$ whose intersection with the horosphere is $B$ and whose projection back to the horosphere is smallest is a horoball, of infinite radius. The radius of its projection is $\frac{1}{2}(R^2+1)$. In order to be mapped to disjoint horoballs, two Euclidean $R$-balls $B_1$ and $B_2$ must have centers at distance at least $R^2 +1$. This makes the scaling factor $\ell$ depend on $R$. Thus this embedding is not large scale conformal.

\textbf{Question}. Are orbital maps of subgroups in nilpotent groups always large scale conformal ?

\subsection{Assouad's embedding}

Assouad's Embedding Theorem, \cite{A} (see also \cite{NN}), states that every snowflake $X^\alpha=(X,d_X^\alpha)$, $0<\alpha<1$, of a doubling metric space admits a bi-Lipschitz embedding into some Euclidean space. 

\begin{prop}
Bi-Lipschitz embeddings of snowflakes are large scale conformal. In particular, the Assouad embedding of a doubling metric space into Euclidean space is large scale conformal.
\end{prop}

\begin{pf}
Identity $X\to X^\alpha$ is large scale conformal. Indeed, the correspondence $B\to B'$ is the identity, concentric ball $\ell B$ is mapped to $\ell^\alpha B'$. So given a scaling factor $\ell'>1$ in the range, $\ell=\ell'^{1/\alpha}$ fits as a scaling factor in the domain.

Bi-Lipschitz maps are large scale conformal, so the composition is large scale conformal as well (Proposition \ref{composition}).
\end{pf}

\subsection{Sphere packings}

Following Benjamini and Schramm \cite{BS}, we view a sphere packing in $\R^d$ as a map from the vertex set of a graph, the incidence graph $G$ of the packing, to $\R^d$. $G$ carries a canonical packing by balls of radius $\frac{1}{2}$ (whose incidence graph is $G$ itself). It is a $(1,\frac{1}{2},\frac{1}{2})$-packing, which is mapped to the given $(1,0,\infty)$-packing. 

Coarse conformality is a very strong restriction on a sphere packing. It has something to do with uniformity in the sense of \cite{BC}. Recall that a sphere packing is $M$-uniform if 
\begin{enumerate}
  \item the degree of the incidence graph is bounded by $M$,
  \item the ratio of radii of adjacent spheres is $\leq M$.
\end{enumerate}
But uniformity is not sufficient. For instance, let $\sigma(z)=az$ be a planar similarity. If $a=re^{i\theta}$ is suitably chosen (for intance, $\theta=\pi/10$ and $r<e^{-2\pi/100}$), there exists a circle $C$ such that $\sigma(C)$ touches $C$ but no iterate $\sigma^k (C)$, $k\not=1$, $-1$, does. The collection of iterates $\sigma^k (C)$, $k\in\Z$, constitutes a uniform planar circle packing whose incidence graph is $\Z$ but which does not give rise to a roughly conformal map of $\Z$ to $\R^2$. Indeed, if $R$ is large enough so that iterates $\sigma^k (C)$, $k\in\{-R,...,R\}$ make a full turn around the origin, any Euclidean ball $B'$ which contains $2R+1$ consecutive circles of the packing contains the origin. So the image of any $(N,1,R,R)$-packing of $\Z$ has infinite multiplicity at the origin.

However, the same construction with $a>0$ gives rise to a coarse conformal map, composition of an inversion with the standard isometric embedding of $\Z$ in $\R^2$.

\medskip

Every subgraph of the incidence graph of a sphere packing is again the incidence graph of a sphere packing. No such property holds for coarsely (resp. roughly, resp. large scale) conformal mappings.

\section{Quasi-symmetry structures}
\label{qsstructure}

This notion is needed in order to host an important example of rough conformal map, the Poincar\'e model of a hyperbolic metric space.

\subsection{Definition}

\begin{defi}
A \emph{quasi-symmetry} structure on a set $X$ is the data of a set $\mathcal{B}$ with a family $(\Phi_\ell)_{\ell\in\R_+^*}$ of maps $\mathcal{B}\to \mathrm{Subsets}(X)$ satisfying
\begin{align*}
\ell\leq \ell' \quad\implies  \quad \Phi_{\ell}(B)\subset\Phi_{\ell'}(B).
\end{align*}
(in the sequel, one will rarely distinguish an element of $\mathcal{B}$ from the corresponding subset $\Phi_1(B)$ of $X$; then $\Phi_\ell(B)$ will be denoted by $\ell B$).
A set equipped with such a structure is called a \emph{q.s. space}. Elements of $\mathcal{B}$ (as well as the corresponding subsets of $X$) are called \emph{balls}. 
\end{defi}

\begin{ex}
Quasi-metric spaces come with a natural q.s. structure, where $\mathcal{B}=X\times(0,\infty)$, and a pair $(x,r)$ is mapped to the ball $B(x,\ell r)$.
\end{ex}

\begin{defi}
Let $X$ be a q.s. space. An \emph{$(N,\ell)$-packing} is a collection of balls $\{B_j\}$ such that the collection of concentric balls $\ell B_j$ has multiplicity $\leq N$, i.e. the collection $\{\ell B_j\}$ can be split into at most $N$ sub-families, each consisting of pairwise disjoint balls.
\end{defi}

\subsection{Coarsely and roughly conformal maps to q.s. spaces}

\begin{defi}
Let $X$ be a metric space and $X'$ a q.s. space. Let $f:X\to X'$ be a map. Say $f$ is \emph{$(R,S)$-coarsely conformal} if there exist a map 
$$B\mapsto B',\quad\mathcal{B}^{X}_{R,S}\to \mathcal{B}^{X'},$$
and for all $\ell'\geq 1$, an $\ell\geq 1$ and an $N'$ such that
\begin{enumerate}
  \item for all $B\in\mathcal{B}^{X}_{R,S}$, $f(B)\subset B'$.
  \item If $\{B_j\}$ is a $(\ell,R,S)$-packing of $X$, then $\{B'_j\}$ is an $(N',\ell')$-packing of $X'$.
\end{enumerate}
We say that $f$ is \emph{coarsely conformal} if there exists $R>0$ such that for all $S\geq R$, $f$ is $(R,S)$-coarsely conformal.

We say that $f$ is \emph{roughly conformal} if there exists $R>0$ such that $f$ is $(R,\infty)$-coarsely conformal.
\end{defi}

It is harder to be roughly conformal than coarsely conformal.

\subsection{Quasi-symmetric maps between q.s. spaces}

\begin{defi}
Let $X$ and $X'$ be q.s. spaces. A map $f:X\to X'$ is \emph{quasi-symmetric} if for all $\ell'\geq 1$, there exist $\ell\geq 1$ and a correspondence between balls 
$$B\mapsto B',\quad\mathcal{B}^{X}\to \mathcal{B}^{X'},$$
such that for all $B\in\mathcal{B}^{X}$,
\begin{enumerate}
  \item $f(B)\subset B'$.
  \item $f^{-1}(\ell'B')\subset\ell B$.
\end{enumerate}
\end{defi}

Quasi-symmetric maps are morphisms in a category whose objects are q.s. spaces. For q.s. structures associated to metrics, isomorphisms in this category coincide with classical quasi-symmetric homeomorphisms.  

\begin{prop}
Let $X$ be a metric space, and $X'$, $X''$ q.s. spaces. If $f:X\to X'$ is $(R,S)$-coarsely conformal and $f':X'\to X''$ is quasi-symmetric, then the composition $f'\circ f$ is $(R,S)$-coarsely conformal.
\end{prop}

\begin{pf} 
Fix $\ell''\geq 1$. Since $f'$ is quasi-symmetric, there exists $\ell'\geq 1$ such that every ball $B'$ in $X'$ is mapped into a ball $B''$ of $X''$ such that $f'^{-1}(\ell''B'')\subset\ell' B'$. If $B'_1$ and $B'_2$ are balls in $X'$ such that $\ell' B'_1\cap\ell' B'_2=\emptyset$, then $f'^{-1}(\ell''B''_1)\cap f'^{-1}(\ell''B''_2)=\emptyset$, hence $\ell''B''_1$ and $\ell'' B''_2$ are disjoint. Thus $f'\circ f$ is $(R,S)$-coarsely conformal.\end{pf}

\section{The Poincar\'e model}
\label{poinc}

\subsection{The 1-dimensional Poincar\'e model}

The following example of q.s. space may be called the half hyperbolic line.

\begin{defi}
Let $\mathbb{D}$ denote the interval $[0,1]$ equipped with the following q.s. structure. $\mathcal{B}$ is the set of closed intervals in $[0,1]$. For an interval $I\subset(0,1]$ there exists a unique pair $(R,t)$ such that $B=[e^{-R-t},e^{R-t}]$. Then $\ell B:=[e^{-\ell R-t},\min\{1,e^{\ell R-t}\}]$. For an interval of the form $I=[0,b]$, $\ell I=I$.
\end{defi}

In other words, the structure is the usual metric space structure of a half real line transported by the exponential function, and extended a bit arbitrarily to a closed interval.

\begin{rem}
Here is a formula for $\ell I$ when $I=[a,b]$, $a>0$.
\begin{align*}
\ell I=[a^{\frac{1+\ell}{2}}b^{\frac{1-\ell}{2}},\min\{1,a^{\frac{1-\ell}{2}}b^{\frac{1+\ell}{2}}\}].
\end{align*}
\end{rem}
Indeed, if $[a,b]=[e^{-R-t},e^{R-t}]$, then $R=\frac{1}{2}\log(b/a)$ and $t=-\frac{1}{2}\log(ab)$. 

\subsection{Warped products}

There is no way to take products of q.s. spaces, so we use auxiliary quasi-metrics.

\begin{defi}
A \emph{q.m.q.s. space} is the data of a q.s. space and a quasi-metric which defines the same balls (but possibly a different $B\mapsto \ell B$ correspondence). 
\end{defi}

\begin{ex}
We keep denoting by $\mathbb{D}$ the q.m.q.s. space $[0,1]$ equipped with its usual metric but the q.s. structure of $\mathbb{D}$.
\end{ex}

\begin{defi}
The \emph{product} of two q.m.q.s. spaces $Z_1$ and $Z_2$ is $Z_1\times Z_2$ equipped with the product quasi-metric
\begin{align*}
d((z_1,z_2),(z'_1,z'_2))=\max\{d(z_1,z'_1),d(z_2,z'_2)\}
\end{align*}
and, if $B=I\times\beta$ is a ball in the product, $\ell B=\ell I\times \ell \beta$. 
\end{defi}

\begin{ex}
Our main example is $\mathbb{D}\times Z$, where $Z$ is a metric space with its metric q.s. structure.
\end{ex}

\subsection{The Poincar\'e model of a hyperbolic metric space}

\begin{prop}
\label{poinc=>rc}
Let $X$ be a geodesic hyperbolic metric space with ideal boundary $\partial X$. Fix an origin $o\in X$ and equip $\partial X$ with the corresponding visual quasi-distance. Assume that there exists a constant $D$ such that for any $x \in X$ there exists a geodesic ray $\gamma$ from the base point $\gamma(0) = o$ and passing near $x$: $d(x, \gamma) < D$. Then there is a rough conformal map $\pi$, called the \emph{Poincar\'e model}, of $X$ to the q.s. space $\mathbb{D}\times\partial X$. Its image is contained in $(0,1]\times\partial X$. If $x\in X$ tends to $z\in\partial X$ in the natural topology of $X\cup\partial X$, then $\pi(x)$ tends to $(0,z)$.
\end{prop}

\begin{pf}
For $x\in X$, pick a geodesic ray $\gamma$ starting from $o$ that passes at distance $<D$ from $x$. Set $\chi(x)=e^{-d(o,x)}$, $\phi(x)=\gamma(+\infty)$ and $\pi(x)=(\chi(x),\phi(x))$. $\pi$ may be discontinuous, but we do not care. 

Let $B=B(x,R)$ be a ball in $X$, and $t=d(o,x)$. Then $\phi(B)$ is contained in a ball of radius $e^{R-t}$, up to a multiplicative constant. Indeed, thanks to V. Shchur's Lemma, \cite{Shchur}, there exists a constant $C$ such that, if $y\in B$ and $t'=d(o,y)$, then $|t-t'|\leq R+C+2D$ and $t+t'+2\log d_o(\phi(x),\phi(y))\leq R+C+2D$. Thus $2t+2\log d_o(\phi(x),\phi(y))\leq 2R+2C+4D$, whence $d_o(\phi(x),\phi(y))\leq e^{R+C+2D-t}$. It turns out that $\chi(B)$ is also a ball of radius $e^{R-t}$, up to a multiplicative constant. Indeed, if $y\in B$ and $t'=d(o,y)$, then $|t-t'|\leq R$, thus $\chi(y)$ belongs to the interval $[1-e^{R-t},1-e^{-R-t}]$, centered at $1-e^{-t}\cosh(R)$, with radius $e^{-t}\sinh(R)$. From now on, we decide to multiply the quasi-metric on $\partial X$ by the constant factor $e^{-C-2D-1}$, i.e. we use the quasi-metric $d'=e^{-C-2D-1}d_o$. In this way, $\pi(B)$ is contained in the ball $B'$ of $[0,1]\times\partial X$ centered at $(1-e^{-t}\cosh(R),\phi(x))$, with radius $e^{-t}\sinh(R)$. 

Conversely, if $(\tau,\eta)\in B'$, then $\eta\in \partial X$ belongs to $B(\phi(x),e^{-t}\sinh(R))$ (which means that $d'(\phi(x),\eta	)=e^{-C-2D-1}d_o(\phi(x),\eta)\leq e^{-t}\sinh(R)$), and $\tau\in[1-e^{R-t},1-e^{-R-t}]$. Let $t'$ be such that $1-e^{-t'}=\tau$. Then $t-R\leq t'\leq t+R$. Let $y$ (resp. $z$) be the point of the geodesic from $o$ to $\eta$ such that $d(o,y)=t'$ (resp. $d(o,z)=t$). Then $d(x,z)\leq 2t+2\log d_o(\phi(x),\eta)+C+2D\leq 2R+3C+4D+2$, and $d(z,y)=|t-t'|\leq R$. Thus $d(x,y)\leq 4R$ provided $R\geq 3C+4D+2$. We conclude that
\begin{align*}
\pi(B)\subset B', \quad \textrm{and} \quad B'\subset\pi(3B).
\end{align*}

Let us show that
\begin{align*}
\pi(3\ell B)\supset \ell B'.
\end{align*}
To avoid confusion, let us denote by $m_\ell$ the $\R^*$ action on intervals. If $B=B(x,R)$ and $\chi(x)=t$, then $\chi(B)=[e^{-R-t},e^{R-t}]$ and $\chi(\ell B)=m_{\ell}(\chi(B))$ by definition of $m_{\ell}$. On the other hand, $B'=\chi(B)\times\beta$ where $\phi(B)\subset\beta=B(\phi(x),e^{-t}\sinh(R))$, 
\begin{align*}
\phi(\ell B)=B(\phi(x),e^{-t}\sinh(\ell R))\supset B(\phi(x),e^{-t}e^{(\ell-1)R}\sinh(R))\supset\ell B,
\end{align*}
provided $R\geq 1$. Therefore 
\begin{align*}
m_{\ell}(B')=m_{\ell}(\chi(\ell B))\times\ell\beta\subset m_{\ell}(\chi(\ell B))\times \phi(\ell B).
\end{align*}
Since this is the ball $(\ell B)'$ corresponding to $\ell B$, it is contained in $\pi(3\ell B)$. 

Since $X$ is hyperbolic, there is a constant $\delta$ such that two geodesics with the same endpoints are contained in the $\delta$-tubular neighborhood of each other. We can assume that $\delta\leq C$. It follows that if $\pi(x_1)=\pi(x_2)$, then $d(x_1,x_2)\leq\delta$. Let $B_1$ and $B_2$ be balls of radii $\geq C$. Let $B'_1$, $B'_2$ be the corresponding balls in $\mathbb{D}\times\partial X$. If $\ell' B'_1\cap\ell' B'_2\not=\emptyset$, then $\pi(3\ell' B_1)\cap \pi(3\ell' B_2)\not=\emptyset$. Thus $d(3\ell' B_1,3\ell' B_2)\leq\delta\leq C$. If $B_1=B(x_1,r_1)$, $3\ell'B_1=B(x_1,3\ell' r_1)$, its $C$-neigborhood is contained in $B(x_1,3\ell' r_1 +C)\subset B(x_1,4\ell' r_1)=4\ell' B_1$. Hence $4\ell' B_1\cap 4\ell' B_2\not=\emptyset$. We may set $\ell=4\ell'$ and conclude that $\pi$ is $(3C+4D+2,\infty,0,\infty)$-coarsely quasi-symmetric, hence roughly conformal. \end{pf}

\section{Energy}
\label{energy}

Classical analysis has made considerable use of the fact that, on $n$-space, the $L^n$ norm of the gradient of functions is a conformal invariant. Benjamini and Schramm observe that only the dimension of the range matters. In fact, some coarse analogue of the $L^p$ norm of the gradient of functions turns out to be natural under coarsely conformal mappings. This works for all $p$ (this fact showed up in \cite{Pa}).

\subsection{Energy and packings}

Let $X$ be a metric space. Recall that an $(\ell,R,S)$-packing is a collection of balls $\{B_j\}$, each with radius between $R$ and $S$, such that the concentric balls $\ell B_j$ are pairwise disjoint. If $X$ is merely a q.s. space, $\ell$-packings make sense.

Here is one more avatar of the definition of Sobolev spaces on metric spaces. For earlier attempts, see the surveys \cite{He}, \cite{GT}.

\begin{defi}
Let $X$, $Y$ be metric spaces. Let $\ell\geq 1$. Let $u:X\to Y$ be a map. Define its $p$-energy at parameters $\ell$, $R$ and $S\geq R$ as follows. 
\begin{align*}
E_{\ell,R,S}^{p}(u)=\sup\{\sum_{j}\mathrm{diameter}(u(B_j))^p \,;\,(\ell,R,S)\textrm{-packings }\{B_j\}\}.
\end{align*}
\end{defi}

\begin{rem}
The definition of $E_{\ell,0,\infty}^{p}$ extends to q.s. spaces, we denote it simply by
\begin{align*}
E_{\ell}^{p}(u)=\sup\{\sum_{j}\mathrm{diameter}(u(B_j))^p \,;\,\ell\textrm{-packings }\{B_j\}\}.
\end{align*}
\end{rem}

Our main source of functions with finite energy are Ahlfors regular metric spaces.

\begin{prop}
\label{lipenergy}
Let $X$ be a $d$-Ahlfors regular metric space and $u:X\to Y$ a $C^{\alpha}$-H\"older continuous function with bounded support. Then $E_{\ell,0,\infty}^{p}(u)$ is finite for all $p\geq d/\alpha$ and $\ell> 1$. In particular, $E_{\ell,R,S}^{d/\alpha}(u)$ is finite for all $R,S$.
\end{prop}

\begin{pf}
Let $\mu$ denote Hausdorff $d$-dimensional measure. Assume first that $\mu(X)<\infty$. For every ball $B$,
\begin{align*}
\mathrm{diam}(u(B))^p &\leq \textrm{const.}\,\mathrm{diam(B)}^{p\alpha}\\
&\leq\textrm{const.}\,\mu(B)^{p\alpha/d}\\
 &\leq \textrm{const.}\,\mu(X)^{-1+p\alpha/d}\mu(B),
\end{align*}
thus
\begin{align*}
\sum_{j}\mathrm{diam}(u(B_j))^p &\leq \textrm{const.}\,\mu(X)^{-1+p\alpha/d}\sum_{j}\mu(B_j)\\
&\leq \textrm{const.}\,\mu(X)^{p\alpha/d}<\infty.
\end{align*}
If $X$ is unbounded, we assume that $u$ has support in a ball $B_0$ of radius $R_0$. Given $\ell>1$, if a ball $B$ intersects $B_0$, and has radius $R>2R_0/(\ell-1)$, then $\ell B$ contains $B_0$. If such a ball arises in an $\ell$-packing, it is the only element of the packing, hence an upper bound on $\sum\mathrm{diameter}(u(B))^p \leq (\max u -\min u)^{p}$. Otherwise, all balls of the $\ell$-packing are contained in $(1+(2R_0/(\ell-1))B_0$, whence the above upper bound is valid with $X$ replaced with that ball.\end{pf}

\begin{ex}
\label{extwisted}
Let $0<\alpha\leq 1$. Let $Z$ be a compact $Q$-Ahlfors regular metric space. Let $Z^\alpha=(Z,d_Z^\alpha)$ be a snowflake of $Z$. Consider the q.m.q.s. space $X=\mathbb{D}\times Z$, with its projections
\begin{align*}
u:X'\to Z,\,u(y,z)=z, \quad\textrm{and}\quad v:X'\to Y,\,v(y,z)=y.
\end{align*} 
Then, for all $\ell$ and for all $p\geq \alpha+Q$,
\begin{align*}
E^{p}_{\ell}(u)<+\infty,\quad E^{p/\alpha}_{\ell}(v)<+\infty.
\end{align*}
\end{ex}

\begin{pf}
Since the previous argument, in the compact case, never uses concentric balls $\ell B$, only the balls $B$ themselves, the difference between $\mathbb{D}\times Z^\alpha$ and the metric space product $[0,1]\times Z^\alpha$ disappears. It is compact and $1+\frac{Q}{\alpha}$-Ahlfors regular with Hausdorff measure $\nu=dt\otimes\mu$. For every ball $B\in \R_+\times Z^\alpha$,
\begin{align*}
\mathrm{diameter}(u(B))=\mathrm{diameter}(B)^{1/\alpha}\sim\nu(B)^{\frac{1}{\alpha}\frac{1}{1+\frac{Q}{\alpha}}}=\nu(B)^{\frac{1}{\alpha+Q}},
\end{align*}
thus $u$ has finite $p$-energy for $p\geq \alpha+Q$.\end{pf}

\subsection{Dependence on radii and scaling factor}
\label{radii}

The following lemma shows that for spaces with bounded geometry, taking the supremum over $(\ell,R,S)$-packings does not play such a big role, provided $R>0$ and $S<\infty$. However, $(\ell,0,\infty)$-packings in the definition of $E^{p}_{\ell,0,\infty}$-energy are much harder to handle.

\begin{defi}
A metric space $X$ has the \emph{tripling property} if
\begin{enumerate}
  \item Balls are connected.
  \item There is a function $N'=N'(N,\ell,R,S)$ (called the \emph{tripling function} of $X$) with the following property. For every $N,\ell\geq 1,R>0,S\geq R$, every $(N,\ell,R,S)$-packing of $X$ is an $(N',\frac{R+S}{R},R,S)$-packing as well.
\end{enumerate}
A $(N,\ell,R,S)$-packing of $X$ is \emph{covering} if the interiors form an open covering of $X$.
\end{defi}

The tripling property is meaningful only for $\ell<\frac{R+S}{R}$.
If a metric space $X$ carries a measure $\mu$ such that the measure of $R$-balls is bounded above and below in terms of $R$ only (this will be called a \emph{metric space with controlled balls} below), then it has the tripling property.

\begin{lem}
\label{onepacking}
Let $X$ be a metric space which has the tripling property. Let $\{B_j\}$ be a covering $(N,\ell,R,S)$-packing of $X$. Then for all maps $u:X\to Y$ to metric spaces,
\begin{align*}
E^{p}_{\ell,R,S}(u)\leq N'(N,\ell,R,S)^p\,\sum_{j}\mathrm{diam}(u(B_j))^p ,
\end{align*}
where $N'$ is the tripling function of $X$.
\end{lem}

\begin{pf}
Let $\{B_j\}$ be a covering $(N,\ell,R,S)$-packings of $X$ and $\{B'_k\}$ a $(\ell,R,S)$-packing of $X$. Let $B'$ be a ball from the second packing. Let $J(B')$ index the balls from $\{B_j\}$ whose interiors intersect $B'$. By assumption, $|J(B')|\leq N'(N,\ell,R,S)$. Indeed, if $B$ intersects $B'$, then the center of $B'$ belongs to $\frac{R+S}{R}B$. Since $B'$ is connected, given any two points $x$, $x'\in B'$, there exists a chain $B_{j_1},\ldots,B_{j_k}$, $j_i \in J(B')$ with $x\in B_{j_1}$, $x'\in B_{j_k}$, and each $B_{j_i}$ intersects $B_{j_{i+1}}$. This implies that $d(u(x),u(x'))\leq\sum_{i=1}^{k}\mathrm{diam}(u(B_{j_i}))$, and
\begin{align*}
\mathrm{diam}(u(B'))\leq\sum_{j\in J(B')}\mathrm{diam}(u(B_j)).
\end{align*}
H\"older's inequality yields
\begin{align*}
\mathrm{diam}(u(B'))^p\leq N'^{p-1}\sum_{j\in J(B')}\mathrm{diam}(u(B_j))^p .
\end{align*}
A given ball of $\{B_j\}$ appears in as many $J(B')$ as its interior intersects balls of $\{B'_k\}$. This happens at most $N'(N,\ell,R,S)$ times. Therefore
\begin{align*}
\sum_{k}\mathrm{diam}(u(B'_k))^p\leq N'^{p}\sum_{j}\mathrm{diam}(u(B_j))^p .
\end{align*}
\end{pf}

\begin{cor}
\label{rad}
Let $X$ be a metric space which has the tripling property. Up to multiplicative constants depending only on $R$, $S$ and $\ell$, energies $E^{p}_{\ell,R,S}(u)$ do not depend on the choices of $R>0$, $S<\infty$ and $\ell\geq 1$.
\end{cor}

\begin{pf}
Let us show that covering $(N,\ell,R,R)$-packings exist for all $R>0$ and $\ell\geq 1$, provided $N$ is large enough, $N=N(\ell,R)$. Let $\{B_j=B(x_j,\frac{R}{2})\}$ be a maximal collection of disjoint $\frac{R}{2}$-balls in $X$. In particular, $\{B_j\}$ is a $(1,\frac{R}{2},(2\ell-1)\frac{R}{2})$-packing. By the tripling property, there exists $N=N'(1,2\ell,\frac{R}{2},(2\ell-1)\frac{R}{2})$ such that $\{B_j\}$ is a $(N,2\ell,\frac{R}{2},\frac{R}{2})$-packing. Then $\{2B_j\}$ is a covering $(N,\ell,R,R)$-packing.

Fix $0<R\leq 1\leq S<+\infty$. Let $\{B_j\}$ be a covering $(N(\ell,1),\ell,1,1)$-packing. According to Lemma \ref{onepacking}, for all maps $u$ to metric spaces,
\begin{align*}
E^{p}_{\ell,1,1}(u)&\leq E^{p}_{\ell,R,S}(u)\\
&\leq N'(1,\ell,R,S)^p \sum_{j}\mathrm{diam}(u(B_j))^p \\
&\leq N'(1,\ell,R,S)^p N(\ell,1) E^{p}_{\ell,1,1}(u),
\end{align*}
so changing radii is harmless. By definition of tripling, given $\ell'\geq\ell$, every $(\ell,1,\ell')$-packing is simultaneously a $(N'',\ell'+1,1,\ell')$-packing, $N''=N'(1,\ell,1,\ell')$, hence
\begin{align*}
E^{p}_{\ell',1,1}(u)&\leq E^{p}_{\ell,1,1}(u)\\
&\leq E^{p}_{\ell,1,\ell'}(u)\\
&\leq N'' E^{p}_{\ell'+1,1,\ell'}(u) \\
&\leq N'' E^{p}_{\ell',1,\ell'}(u)\\
&\leq N'' N'(1,\ell',1,\ell')^p N(\ell',1)E^{p}_{\ell',1,1}(u),
\end{align*}
so changing scaling factor is also harmless.
\end{pf}

\medskip

If $p=1$, even the upper bound on radii of balls does not play such a big role.

\begin{lem}
Let $X$ be a geodesic metric space. For every function $u$ on $X$,
\begin{align*}
E^{1}_{\ell,R,R}(u)\leq E^{1}_{\ell,R,\infty}(u)\leq (2\ell+2) E^{1}_{\ell,R,R}(u).
\end{align*}
\end{lem}

\begin{pf}
Let $B$ be a large ball. Assume $u$ achieves its maximum on $B$ at $x$ and its minimum at $y$. Along the geodesic from $x$ to $y$, consider an array of touching $R$-balls $B_j$. Then
\begin{align*}
\mathrm{diameter}(u(B))\leq\sum_{j}\mathrm{diameter}(u(B_j)).
\end{align*} 
Pick one ball every $2\ell$ along the array, in order to get an $(\ell,R,R)$-packing. The array is the union of $2\ell$ such packings, whence 
$$\sum_{j}\mathrm{diameter}(u(B_j))\leq E^{1}_{\ell,R,R}(u).$$

Summing up over balls of an arbitrary $(\ell,R,\infty)$-packing yields the announced inequality.
\end{pf}

\subsection{Functoriality of energy}

\begin{lem}
\label{energytransport}
Let $X$, $X'$ and $Y$ be metric spaces. Let $f:X\to X'$ be $(R,S,R',S')$-coarsely conformal. Then for all $\ell'\geq 1$, there exists $\ell\geq 1$ and $N'$ such that for all maps $u:X'\to Y$, 
\begin{align*}
E_{\ell, R,S}^{p}(u\circ f)\leq N'\,E_{\ell',R',S'}^{p}(u).
\end{align*}
\end{lem}

\begin{pf}
Let $\{B_j\}$ be an $(\ell,R,S)$-packing of $X$. Let $\{B'_j\}$ be the corresponding $(N',\ell',R',S')$-packing of $X'$. Split it into $N'$ sub-collections which are $(1,\ell',R',S')$-packings. By assumption, $f(B_j) \subset B'_j$, so $\mathrm{diam}(u\circ f(B_j))\leq \mathrm{diam}(u(B'_j))$. This yields, for each sub-collection,
\begin{align*}
\sum_{j}diam(u\circ f(B_j))^p \leq\sum_{j}diam(u(B'_j))^p \leq E_{\ell',R',S'}^{p}(u).
\end{align*}
Summing up and taking supremum, this shows that $E_{\ell, R,S}^{p}(u\circ f)\leq N'\,E_{\ell',R',S'}^{p}(u)$.\end{pf}

\begin{ex}
If $Y$ is $d$-Ahlfors regular and compact, then the identity $Y\to Y$ has finite $E_{\ell', R',S'}^{p}$-energy for all $\ell', R',S'$, so coarsely conformal maps $X\to Y$ have finite $E_{\ell, R,S}^{p}$-energy themselves for suitable $\ell$. 
\end{ex}

\subsection{\texorpdfstring{$(1,1)$}{}-curves}

\begin{defi}
A \emph{$(1,1)$-curve} in a metric space $X$, is a map $\gamma:\N\to X$ which is $(1,1,R,S)$-coarsely conformal for all $R>0$ and $S\geq R$. 
When $X$ is locally compact and equipped with a base point $o$, a \emph{based $(1,1)$-curve} is a proper $(1,1)$-curve such that $\gamma(0)=o$. 
\end{defi}

An $(\ell,1,1)$-packing of $\N$ corresponds to a subset $A\subset\N$ such that every $\ell$-ball centered at a point of $A$ contains at most one points of $A$. The packing consists of unit balls centered at points of $A$. Let us call such a set an \emph{$\ell$-subset} of $\N$.
A $(1,1)$-curve in $X$ is a sequence $(x_i)_{i\in\N}$ such that for all $R$ and all sufficiently large $S\gg R$, there exists a collection of balls $B_i$ in $X$ with radii between $R$ and $S$ such that 
\begin{itemize}
  \item $B_i$ contains $\{x_{i-1},x_{i},x_{i+1}\}$.
  \item For all $\ell'\geq 1$, there exist $\ell$ ans $N'$ such that for every $\ell$-subset of $\N$, the collection of concentric balls $\{\ell'B'_i\,;\,i\in A\}$ has multiplicity $\leq N'$.
\end{itemize}
Thus a $(1,1)$-curve is a chain of slightly overlapping balls which do not overlap too much: if radii are enlarged $\ell'$ times, decimating (i.e. keeping only one ball every $\ell$) keeps the collection disjoint or at least bounded multiplicity.

\begin{ex}
An isometric map $\N\to X$ is a $(1,1)$-curve. In particular, geodesic rays in Riemannian manifolds give rise to $(1,1)$-curves.
\end{ex}
This is a special case of Lemma \ref{coarse=>lsc}.

\begin{rem}
In the definition of based $(1,1)$-curves, the properness assumption is needed only if $R=0$.
\end{rem}

\begin{pf}
Let $\gamma:\N\to X$ be a based $(1,1)$-curve. Let $\ell'=2$, let $\ell$ be the corresponding scaling factor in the domain. The covering of $\N$ by $1$-balls $B(j,1)$ is mapped to balls $B'_j$ of radii $\geq R>0$, with $\gamma(j)\in B'_{j}$. Since the covering $\{B(j,1)\}$ is the union of exactly $\ell$ $\ell$-packings, $\{B'_j\}$ is the union of $\ell$ $(N',2)$-packings, it is an $(\ell N',2)$-packing. In particular, no point of $X$ is contained in more than $\ell N'$ balls $2B'_j$.

If $\gamma$ is not proper, there exists a sequence $i_j$ tending to infinity such that $\gamma(i_j)$ has a limit $x\in X$. Since $B'_{i_j}$ has radius $\geq R$, for $j$ large enough $2B'_{i_j}$ contains $x$, contradicting multiplicity $\leq \ell N'$.\end{pf}

\begin{defi}
Let $X$ be a q.s. space, and $K\subset X$. A \emph{coarse curve} in $X$ is a coarse conformal map $\N\to X$. A coarse curve $\gamma$ is \emph{based at} $K$ if $\gamma(0)\in K$.
\end{defi}

\begin{ex}
Given an isometric map $\gamma:[0,1]\to X$, set $\gamma'(j)=\gamma(\frac{1}{j+1})$. This is a coarse curve in $X\setminus\{\gamma(0)\}$. In particular, geodesic segments in punctured Riemannian manifolds give rise to $(0,\infty)$-curves.
\end{ex}
Since $\gamma'$ is the composition of inversion $\R\to\R$ and an isometric embedding, this is a special case of Proposition \ref{qs=>cc}.

\begin{prop}
Let $X$ be a metric space and $X'$ a q.s. space. If $\gamma:\N\to X$ is a $(1,1)$-curve and $f:X\to X'$ is a coarse conformal map, then $f\circ\gamma$ is a coarse curve.
\end{prop}

\subsection{Modulus}

We need to show that certain maps with finite energy have a limit along at least one based curve. To do this, we shall use the idea, that arouse in complex analysis, of a property satisfied by almost every curve.

\begin{defi}
Let $Y$ be a metric space. The \emph{length} of a map $u:\N\to Y$ is
\begin{align*}
\mathrm{length}(u)=\sum_{i=0}^{\infty}d(u(i),u(i+1)).
\end{align*}
\end{defi}

\begin{defi}
\label{defmod}
Let $X$ be a metric space. Let $\Gamma$ be a family of $(1,1)$-curves in $X$. The {\em $(p,\ell,R,S)$-modulus} $\mathrm{mod}_{p,\ell,R,S}(\Gamma)$ is the infimum of $E^{p}_{\ell,R,S}$-energies of maps $u : X\to Y$ to metric spaces such that for every curve $\gamma\in\Gamma$, $\mathrm{length}(u\circ\gamma)\geq 1$.
\end{defi}

\begin{rem}
The definition of $(p,\ell,0,\infty)$-modulus extends to q.s. spaces $X$,
\begin{align*}
\mathrm{mod}_{p,\ell}(\Gamma)=\inf\{E^{p}_{\ell}(u)\,;\,u:X\to Y,\,\mathrm{length}(u\circ\gamma)\geq 1\,\forall\gamma\in\Gamma\}.
\end{align*}
\end{rem}

\begin{lem}
\label{addit}
Let $X$ be a metric space. The union of a countable collection of $(1,1)$-curve families which have vanishing $(p,\ell,R,S)$-modulus also has vanishing $(p,\ell,R,S)$-modulus.
\end{lem}

\begin{pf}
Fix $\epsilon>0$. Let $u_j :X\to \R$ be a function with $E^{p}_{\ell,R,S}(u_j) \leq 2^{-j}\epsilon$ such that for all curves $\gamma$ in the $j$-th family $\Gamma_j$, $\mathrm{length}(u_j \circ\gamma)\geq 1$. Consider the $\ell^p$ direct product $Y=\prod_{j}Y_{j}$, i.e.
\begin{align*}
d^{Y}((y_j),(y'_j))=\left(\sum_{j} d(y_j ,y'_j)^p\right)^{1/p},
\end{align*}
and the product map $u=(u_j):X\to Y$. Then 
$$E^{p}_{\ell,R,S}(u)\leq\sum_{j}E^{p}_{\ell,R,S}(u_j)\leq const.\epsilon^p,$$
whereas for all curves $\gamma$ in the union curve family, 
\begin{align*}
\mathrm{length}(u\circ\gamma)\geq \sup_j \mathrm{length}(u_j \circ\gamma)\geq 1.
\end{align*}
This shows that the union curve family has vanishing $(p,\ell,R,S)$-modulus.\end{pf}

\begin{lem}
\label{mod0}
Let $X$ be a metric space. Let $\Gamma$ be a family of $(1,1)$-curves in $X$. Then $\Gamma$ has vanishing $(p,\ell,R,S)$-modulus if and only if there exists a function $u:X\to\R_+$ such that $E^{p}_{\ell,R,S}(u)<+\infty$ but $\mathrm{length}(u\circ\gamma)=+\infty$ for every $\gamma\in\Gamma$.
\end{lem}

\begin{pf}
One direction is obvious. In the opposite direction, first observe that by rescaling target metric spaces, one can assume that there exist maps $u_j :X\to Y_j$ such that for all $\gamma\in\Gamma$, $\mathrm{length}(u_j\circ\gamma)\geq j$ and $E^{p}_{\ell,R,S}(u_j)<2^{-j}$. Apply the $\ell^p$-product construction again. Get $u=(u_j):X\to Y$ such that $E^{p}_{\ell,R,S}(u)\leq 1$ and $\mathrm{length}(u\circ\gamma)\geq \max_j \mathrm{length}(u_j \circ\gamma)=+\infty$.\end{pf}

\begin{lem}
\label{modlimitf}
Let $X$ be a metric space. Let $Y$ be a complete metric space. Fix $R\leq S$ and $\ell\geq 1$. Let $u:X\to Y$ be a map of finite $E^{p}_{\ell,R,S}$ energy. The family of $(1,1)$-curves along which $u$ does not have a limit has vanishing $(p,\ell,R,S)$-modulus.
\end{lem}

\begin{pf}
If $\mathrm{length}(u\circ\gamma)<\infty$, then $u\circ\gamma$ has a limit in $Y$, since $Y$ is complete. Let $\Gamma_{nl}$ be the sub-family of curves along which the length of $u$ is infinite. It contains all curves along which $u$ does not have a limit. By assumption, $E^{p}_{\ell,R,S}(u)<\infty$, but $\mathrm{length}(u \circ\gamma)=+\infty\geq 1$ for all curves $\gamma\in \Gamma_{nl}$. So $\mathrm{mod}_{p,\ell,R,S}(\Gamma_{nl})=0$. \end{pf}

\begin{rem}
The cases $R=0,S=\infty$ of Lemmata \ref{addit} to \ref{modlimitf} extend to q.s. spaces $X$.
\end{rem}

\subsection{Parabolicity}

\begin{defi}
Let $X$ be a locally compact metric space. Say $X$ is \emph{$(p,\ell,R,S)$-parabolic} if the family of all $(1,1)$-curves based at some point has vanishing $(p,\ell,R,S)$-modulus. 

If $X$ is merely a locally compact q.s. space, $(p,\ell)$-parabolicity means that the $(p,\ell)$-modulus of the family of proper coarse curves based at any compact set with nonempty interior vanishes.
\end{defi}

\begin{rem}
If metric space $X$ has the tripling property, $(p,\ell,R,S)$-parabo\-licity does not depend on the choices of $R>0$, $S<\infty$ and $\ell\geq 1$, thanks to Corollary \ref{rad}.
\end{rem}
On the other hand, $(p,\ell)$-parabolicity may depend wether $\ell=1$ or $\ell>1$, as we shall see in the next subsections.

\begin{prop}
Let $X$ be a metric space and $X'$ a q.s. space. Let $f:X\to X'$ be a coarse conformal map. Let $R$ be large enough, and $S\geq R$. Let $\Gamma$ be a family of $(1,1)$-based curves in $X$. Then for all $\ell'$, there exist $\ell$ such that
\begin{align*}
\mathrm{mod}_{p,\ell,R,S}(\Gamma)\leq \mathrm{mod}_{p,\ell}(f(\Gamma)).
\end{align*} 
\end{prop}

\begin{pf}
By the composition rule (Lemma \ref{composition}), $f(\Gamma)$ is a family of $(0,\infty)$-based curves. Given $u:X'\to Y$ such that $\mathrm{length}(u\circ \gamma')\geq 1$ for all $\gamma'\in f(\Gamma)$, $\mathrm{length}(u\circ f\circ\gamma)\geq 1$ for all $\gamma\in \Gamma$. According to Lemma \ref{energytransport}, for all $\ell'\geq 1$, there exists $\ell\geq 1$ such that
\begin{align*}
E_{\ell, R,S}^{p}(u\circ f)\leq E_{\ell'}^{p}(u).
\end{align*}
Taking the infimum over such maps $u$ yields the announced inequality. \end{pf}

\begin{cor}
\label{propercc=>para}
Let $X$ be a metric space and $X'$ a q.s. space. Let $f:X\to X'$ be a proper, coarse conformal map. Then there exists $R>0$ such that for all $S\geq R$, if $X'$ is $(p,\ell')$-parabolic for some $\ell'$, then $X$ is $(p,\ell,R,S)$-parabolic for a suitable $\ell$.
\end{cor}

There is a similar statement for $(R',S',R,S)$-coarsely conformal maps. This shows that parabolicity does not depend on the choice of base point, provided one accepts to change parameters $\ell$, $R$ and $S$. Indeed, the map which is identity but for one point $o$ which is mapped to $o'$ is proper and $(R',S',R,S)$-coarsely conformal.

\begin{cor}
\label{lsc=>para}
Let $X$ and $X'$ be a metric spaces. Let $f:X\to X'$ be a uniformly conformal map. For every $R'>0$, there exists $R$ such that for all $S\geq R$, if $X'$ is $(p,\ell',R',\infty)$-parabolic for some $\ell'$, then $X$ is $(p,\ell,R,S)$-parabolic for a suitable $\ell$.
\end{cor}

\begin{pf}
Uniformly conformal maps are proper. 
\end{pf}

\subsection{Parabolicity of Ahlfors-regular spaces}

\begin{prop}
\label{parabolicnc}
Let $X$ be a non-compact $Q$-Ahlfors regular metric space with $Q>1$. Let $K$ be a ball. For all $\ell>1$, there exists a finite $(Q,\ell)$-energy function $w:X\to\R$ which has no limit along every coarse curve based on $K$. It follows that $X$ is $(p,\ell)$-parabolic for every $\ell>1$ and every $p\geq Q$.
\end{prop}

\begin{pf}
Let $\mu$ be a measure such that balls of radius $\rho$ have measure $\rho^Q$ up to multiplicative constants.

Let $m=\max\{\frac{\ell+1}{\ell-1},e\}$. Fix an origin $o\in X$, set $r(x)=d(x,o)$, $v(r)=\log\log r$ and $w(x)=\sin(v(r(x)))$ if $r\geq m^2$, $=\log\log(m^2)$ otherwise. Let $\{B_j\}$ be a $\ell$-packing of $X$. At most one ball $B$ is such that $o\in\ell B$, it contributes to $\sum_j \mathrm{diameter}(w(B_j))^p$ by at most 1. We shall ignore it henceforth. Other balls $B=B(x,\rho)$ are such that $o\notin\ell B$, hence $r(x)=d(o,x)>\ell \rho$. For all $x'\in B$, $r(x')=d(o,x')$ satisfies $r(x)-\rho\leq r(x')\leq r(x)+\rho$, hence
\begin{align*}
\frac{\sup_B r}{\inf_B r}\leq\frac{r(x)+\rho}{r(x)-\rho}=\frac{\frac{r(x)}{\rho}+1}{\frac{r(x)}{\rho}-1}\leq\frac{\ell+1}{\ell-1}\leq m.
\end{align*}
For $i\in\Z$, let $r_i=m^{i}$ and define $Y_i=\{x\in X\,;\,r_{i-2}\leq r(x)\leq r_{i}\}$ and $L_i=Lip(v_{|Y_i})$. Note that $\mu(Y_i)\leq C\,r_i^p$. By construction, each ball $B$ of the packing is contained in at least of the $Y_i$. If $i\leq 2$, $w$ is constant on $Y_i$, such balls do not contribute. Let $i\geq 3$. For a ball $B\subset Y_i$,
\begin{align*}
\mathrm{diameter}(w(B))^Q&\leq&\mathrm{diameter}(v(B))^Q\\
&\leq&L_i^Q \mathrm{diameter}(r(B))^Q\\
&\leq&C\,L_i^Q \mu(B).
\end{align*}
Summing up over all balls of the packing contained in $Y_i$,
\begin{align*}
\sum_{B_j\subset Y_i}\mathrm{diameter}(w(B_j))^Q \leq C\,L_i^Q \mu(Y_i)\leq C\,(L_i r_i)^Q.
\end{align*}
Since 
$v'(t)=\frac{1}{t\log t}$, $L_i\leq\frac{1}{r_{i-2}\log(r_{i-2})}$, $L_i r_i\leq \frac{r_i}{r_{i-2}\log(r_{i-2})}=\frac{m^2}{(i-2)\log m}$. Summing up over $i\geq 3$ leads to
\begin{align*}
\sum_{j}\mathrm{diameter}(B_j)^Q\leq C\,\sum_{i=3}^{\infty}(\frac{1}{i-2})^Q<\infty.
\end{align*}
This shows that $E^{p}_{\ell}(w)<\infty$.

Let $\gamma:\N\to X$ be a proper coarse curve based at $K=B(o,\rho)$. Let us show that $w\circ\gamma$ has no limit. By definition, there exists $\tilde{\ell}\geq 1$ such that every $\tilde{\ell}$-packing of $\N$ is mapped to a $\ell$-packing $\{B'_j\}$ of $X$. For each $i=0,\ldots,\ell-1$, this applies to the $\ell$-packing of unit balls centered at $3\tilde{\ell}\N+i\subset\N$. We know that
\begin{align*}
\sum_{j\in\N}\mathrm{diameter}(v(B'_{3\tilde{\ell} j+i}))^{p}<+\infty,
\end{align*}
hence, summing over $i$,
\begin{align*}
\sum_{j\in\N}\mathrm{diameter}(v(B'_{j}))^{p}<+\infty,
\end{align*}
so these diameters tend to zero. Since $B'_j$ contains $\gamma(j)$ and $\gamma(j+1)$, this implies that $|v\circ\gamma(j+1)-v\circ\gamma(j)|$ tends to 0. Therefore the $\omega$-limit set of the sequence $w\circ\gamma$ is the whole interval $[-1,1]$, $w\circ\gamma$ has no limit. 

We conclude that the family of all proper coarse curves based at some ball has vanishing $(p,\ell)$-modulus, i.e., $X$ is $(p,\ell)$-parabolic. A fortiori, the family of all proper coarse curves based at $o$ has vanishing $(p,\ell,R,S)$-modulus, i.e., $X$ is $(p,\ell,R,S)$-parabolic, for $R>0$.
\end{pf}

\begin{rem}
\label{lsparabolic}
If we are merely interested in $(p,\ell,R,S)$-parabolicity for some $R>0$, or $(p,\ell,R,\infty)$-parabolicity, a weaker form of Ahlfors-regularity is required. It suffices that balls of radius $\rho\geq R$ satisfy $c\,\rho^{Q}\leq \mu(B)\leq C\,\rho^Q$. Let us call this \emph{$Q$-Ahlfors regularity in the large}.
\end{rem}
Indeed, the argument uses only balls of radius $\geq R$.

\begin{cor}
\label{parabolicinthe large}
Let $X$ be a non-compact metric space. Let $Q>1$. Assume that $X$ is $Q$-Ahlfors regular in the large. Then $X$ is $(p,\ell,R,\infty)$-parabolic for every $\ell>1$ and every $p\geq Q$. \emph{A fortiori}, it is $(p,\ell,R,S)$-parabolic for every $\ell>1$, every $0<R\leq S$ and every $p\geq Q$.
\end{cor}

\begin{prop}
\label{parabolic}
Let $X$ be a compact $p$-Ahlfors regular metric space, with $p>1$. Let $x_0\in X$. For every $\ell>1$, there exists a function $w:X\setminus\{x_0\}\to\R$ such that $E_{\ell}^{p}(w)<+\infty$ and $w$ has a limit along no coarse curve converging to $x$.
It follows that $X\setminus\{x_0\}$ is $(p,\ell)$-parabolic.
\end{prop}

\begin{pf}
The same as for the non-compact case, replacing function $r$ with $1/r$.
\end{pf}

\subsection{Parabolicity of \texorpdfstring{$\mathbb{D}$}{}}

The half real line is 1-Ahlfors regular, a case which is not covered by Proposition \ref{parabolicnc}. It is not 1-parabolic. Indeed, any function of finite 1-energy on $\R_+$ has a limit at infinity. However, it is $p$-parabolic for every $p>1$.

\begin{lem}\label{rparabolic}
The half real line $\R_+$ equipped with its metric q.s. structure is $p$-parabolic for every $p>1$.
\end{lem}

\begin{pf}
Let $\ell>1$. Denote by $m=\frac{\ell+1}{\ell-1}$. We can assume that $m\geq e$. Define 
\begin{align*}
u(t)=\log\log |t| \textrm{ if }t\geq m, \quad
u(t)=\log\log m   & \textrm{otherwise}.
\end{align*}
Let us show that $u$ has finite $(p,\ell)$-energy for all $p>1$ and $\ell>1$. Let $\{B_j\}$ be an $\ell$-packing of $\R_+$. By assumption, the collection of concentric balls $\{\ell B_j\}$ consists of disjoint intervals. For simplicity, assume that $0$ belongs to one of the $\ell B_j$'s, say $\ell B_0$ (otherwise, translate everything). Write $B_j=[a_j,b_j]$ and assume that $a_j\geq 0$. Since $\ell B_0$ and $\ell B_j$ are disjoint,
\begin{align*}
\frac{a_j+b_j}{2}-\ell\frac{b_j-a_j}{2}\geq 0,
\end{align*}
thus
\begin{align*}
b_j\leq\frac{\ell+1}{\ell-1}a_j =m a_j.
\end{align*}

We split the sum $\sum_j |u(b_j)-u(a_j)|^p$ into sub-sums where $a_j\in[m^i,m^{i+1})$. Since intervals $B_j$ are disjoint and $u$ is nondecreasing,
\begin{align*}
\sum_{a_j\in[m^i,m^{i+1})}|u(b_j)-u(a_j)|&\leq u(m^{i+2})-u(m^i)\\
&= \log((i+2)\log m)-\log(i\log m)\\
&= \log\frac{i+2}{i}\leq\frac{2}{i}.
\end{align*}
Next, we use the general inequality, for nonnegative numbers $x_j$,
\begin{align*}
\sum x_j^p\leq(\sum x_j)^p,
\end{align*}
and get
\begin{align*}
\sum_{a_j\in[m^i,m^{i+1})}|u(b_j)-u(a_j)|^p\leq(\sum_{a_j\in[m^i,m^{i+1})}|u(b_j)-u(a_j)|)^p\leq(\frac{2}{i})^p.
\end{align*}
This gives
\begin{align*}
\sum_{a_j\geq m}|u(b_j)-u(a_j)|^p\leq\sum_{i=1}^{\infty}(\frac{2}{i})^p<+\infty.
\end{align*}
On the remaining intervals, $u$ is constant, except possibly on one inerval containing $m$. On this interval, $|u(b_j)-u(a_j)|\leq u(m^2)$, so its contribution is bounded independently of the packing. We conclude that the supremum over $\ell$-packings of $\sum_{j}|u(b_j)-u(a_j)|^p$ is bounded, i.e. $u$ has finite $p$-energy.\end{pf}

Since, as a q.s. space, $\mathbb{D}$ is isomorphic to the half real line, $\mathbb{D}$ is $p$-parabolic for all $p>1$ as well. As a preparation for the next result, note that the isomorphism is the exponential map $t\mapsto \exp(-t):\R_+\to\mathbb{D}$. Therefore the function of finite energy on $\mathbb{D}$ is $w(y)=\sin\log|\log|\log y||$.

\subsection{Parabolicity of warped products}

\begin{prop}
\label{twisted}
Let $Z$ be a compact $Q$-Ahlfors regular metric space. Let $0<\alpha\leq 1$. Let $z_0\in Z$. Then $\mathbb{D}\times Z^{\alpha}\setminus\{(0,z_0)\}$ is $(\alpha+Q)$-parabolic.
\end{prop}

\begin{pf} 
For $(y,z)\in[0,1]\times Z$, let $r(y,z)=\min\{y,d^{Z^\alpha}(z,z_0)\}$ denote the distance to $(0,z_0)$. We use the bounded function $w=\sin v(r)$, where
\begin{align*}
v(t)=\log|\log|\log t|| \quad \textrm{if}\quad t\leq r_1,\quad v(t)=\log|\log|\log r_1||\quad \mathrm{otherwise}.
\end{align*}
The constant $r_1=r_1(\ell)$ is produced by the following Lemma.

\begin{lem}
\label{r1}
Let $v(t)=\log|\log|\log t||$. Assume that $\ell>1$, and that
\begin{align*}
b\leq r_1:=\min\{\frac{\ell-1}{\ell},\ell^{-2},e^{-e^2}\}.
\end{align*}
Then, if $\displaystyle\frac{b}{a}\geq\frac{\ell}{\ell-1}$,
\begin{align*}
v(a)-v(\ell b)\leq 16\,(v(a)-v(b)).
\end{align*}
\end{lem}

\begin{pf}
We use the inequalities
\begin{align*}
0\leq u\leq 1 \implies \log(1+u)\geq \frac{1}{2}u,\quad
0\leq u\leq \frac{1}{2} \implies -\log(1-u)\leq 2u.
\end{align*}
Set $t=\log\frac{1}{b}$ and $s'=\log\frac{1}{a}-\log\frac{1}{b}=\log\frac{b}{a}\geq s:=\log\frac{\ell}{\ell-1}$. Then
\begin{align*}
v(a)-v(b)&= \log\log\log\frac{1}{a}-\log\log\log\frac{1}{b}\\
&= \log\log(t+s)-\log\log(t)\\
&= \log\left(\frac{\log(t+s')}{\log(t)}\right)\\
&\geq \log\left(\frac{\log(t+s)}{\log(t)}\right)\\
&= \log\left(1+\frac{\log(t+s)-\log(t)}{\log(t)}\right)\\
&= \log\left(1+\frac{\log(1+(s/t))}{\log(t)}\right)\\
\end{align*}
If $t=\log\frac{1}{b}\geq\max\{s,e^2\}$, $s/t\leq 1$, $\log(t)\geq 2$, so $\frac{\log(1+(s/t))}{\log(t)}\leq 1$. Also $s/t\leq 1$, thus
\begin{align*}
\log(1+(s/t))\leq \frac{1}{2}\frac{s}{t},\quad\log\left(1+\frac{\log(1+(s/t))}{\log(t)}\right)\leq\frac{1}{2}\frac{\log(1+(s/t))}{\log(t)}\leq\frac{1}{4}\frac{s}{t\log(t)},
\end{align*}
and
\begin{align*}
v(a)-v(b)\geq \frac{s}{4t\log(t)}.
\end{align*}
Conversely, setting $\sigma=\log\frac{1}{b}-\log\frac{1}{\ell b}=\log\ell$,
\begin{align*}
v(b)-v(\ell b)&= \log\log\log\frac{1}{b}-\log\log\log\frac{1}{\ell b}\\
&= \log\log(t)-\log\log(t-\sigma)\\
&= -\log\left(\frac{\log(t-\sigma)}{\log(t)}\right)\\
&= -\log\left(1+\frac{\log(t-\sigma)-\log(t)}{\log(t)}\right)\\
&= -\log\left(1+\frac{\log(1-(\sigma/t))}{\log(t)}\right)\\
&\leq -2\frac{\log(1-(\sigma/t))}{\log(t)}\\
&\leq 4\frac{\sigma}{t\log(t)},
\end{align*}
provided $\sigma/t\leq \frac{1}{2}$ and $-\frac{\log(1-(\sigma/t))}{\log(t)}\leq \frac{1}{2}$, which holds if $t=\log\frac{1}{b}\geq\max\{2\sigma,e^2\}$.
\end{pf}

\bigskip

\textbf{Back to the proof of Proposition \ref{twisted}}.

Let $\{B_j\}$ be a $\ell$-packing of $\mathbb{D}\times Z$.
The packing splits into three sub-collections, 
\begin{enumerate}
  \item Balls that contain $(0,z_0)$. 
  \item Balls intersecting $\R_+\times \{z_0\}$, but not containing $(0,z_0)$.
  \item Balls which do not intersect $\R_+\times \{z_0\}$.
\end{enumerate}

We first assume that $\alpha=1$, i.e. we ignore any special feature of snowflaked spaces.

\bigskip

1. The first sub-collection has at most one element. Since $|w|\leq 1$, it contributes at most 1 to energy.

\bigskip

2. The second sub-collection is nearly taken care of by Lemma \ref{rparabolic}. 
If two balls $B=I\times\beta$ and $B'=I'\times\beta'$ in $[0,1]\times Z$ are disjoint and both intersect $[0,1]\times\{z_0\}$, then the intervals $I$ and $I'$ are disjoint. In particular, if $\ell B\cap\ell B'=\emptyset$, then $m_\ell(I)\cap m_\ell(I')=\emptyset$. Furthermore, if $B=I\times\beta$ and $I=[a,b]$, the radius of $B$ is $\frac{b-a}{2}$, thus for $z\in\beta$, $d(z,z_0)\leq 2\ell\frac{b-a}{2}$, so
\begin{align*}
\min_B r\geq a,\quad \max_B d(\cdot,z_0)\leq \ell(b-a)\leq \ell b,
\end{align*}
whence
\begin{align*}
\max_B r\leq \ell b.
\end{align*}

The $Z$ factor plays little role. We expect from Lemma \ref{rparabolic} that the sum of $p$-th powers of diameters of images $w(B_j)$ should be bounded independently of the packing, as soon as $p>1$. 

Here comes the proof. First, the last estimate needs be sharpened. Either $\ell(b-a)\leq b$, in which case $r(B)=[a,b]$ and $v((r(B)))=[v(b),v(a)]$, or $\frac{b}{a}\geq\frac{\ell}{\ell-1}$. Lemma \ref{r1} shows that, in both cases,
\begin{align*}
\mathrm{diameter}(v(r(B)))\leq C\,(v(a)-v(b)),
\end{align*}
provided $b$ is small enough.

\bigskip

If an interval $I=[a,b]$ of $(0,1]$ is such that $m_\ell(I)$ does not contain 1, then $a\geq b^{m}$, for $m=\frac{\ell+1}{\ell-1}$. Indeed, the upper bound of $m_\ell(I)$ is 
$a^{\frac{1-\ell}{2}}b^{\frac{1+\ell}{2}}$. Hence if $B=I\times\beta$ is a ball of $[0,1]\times Z$,
\begin{align*}
\min_B r=a\geq b^{m},\quad \max_B r\leq \ell b\leq \ell(\min_B r)^{1/m}.
\end{align*}
Define inductively a sequence $r_i$ by $r_{i+1}=(\frac{r_i}{\ell})^{m}$ (this gives $\displaystyle r_i=\ell^{-m\frac{m^i-1}{m-1}}$). Also, define $r_0$ so that $r_{1}=(\frac{r_0}{\ell})^{m}$.
Let $\{B_j=I_j\times\beta_j\}$ be a $\ell$-packing of $\mathbb{D}\times Z$. Assume that all $B_j$ intersect $[0,1]\times\{z_0\}$. At most one $\ell B_j$ contains $(1,z_0)$, let us put it aside (it contributes at most 1 to energy). All other intervals $I_j$ are disjoint and each one is contained in at least one of the intervals $[r_{i+2},r_i]$. Therefore the index set is contained in the union of subsets
\begin{align*}
J_i=\{j\,;\, I_j\subset [r_{i+2},r_i]\},
\end{align*} 
and for all $p$,
\begin{align*}
\sum_{j}\mathrm{diameter}(w(B_j))^p \leq \sum_{i=0}^{\infty}\sum_{j\in J_i}\mathrm{diameter}(v(r(B_j)))^p.
\end{align*}

We split the sub-packing into two sub-families,
\begin{enumerate}
  \item Balls $B=I\times\beta$ such that $v(r(B))\not=v(I)$ and $\max I\geq r_1$.
  \item Balls $B=I\times\beta$ such that $v(r(B))=v(I)$ or $\max I\leq r_1$.
\end{enumerate}
The first sub-family satisfy $\max I\geq\frac{\ell}{\ell-1}\min I$ and $\max I\geq r_1$. Those which have $\max I\geq r_0$ have $\min I\geq r_1$, $v$ is constant on such balls. The others are contained in $Y_1$. The number of such disjoint intervals is bounded in terms of $\ell$ only. The diameters $\mathrm{diameter}(v(r(B)))\leq v(r_2)$ are bounded in terms of $\ell$. Therefore the sum of $p$-th powers $\mathrm{diameter}(v(r(B)))$ over this sub-family is apriori bounded in terms of $\ell$ and $p$ only.

\bigskip

The second sub-family satisfy
\begin{align*}
\mathrm{diameter}(v(r(B)))\leq C\,\mathrm{diameter}(v(I)).
\end{align*}
Since the intervals $I_j$ are disjoint,
\begin{align*}
\sum_{j\in J_i}\mathrm{diameter}(w(B_j))^{p}
&\leq \left(\sum_{j\in J_i}\mathrm{diameter}(v(r(B_j)))\right)^{p}\\
&\leq C\,|v(r_{i+2})-v(r_i)|^{p},
\end{align*}
\begin{align*}
\sum_{j}\mathrm{diameter}(w(B_j))^{p}\leq C\,\sum_{i=0}^{\infty}|v(r_{i+2})-v(r_i)|^{p}.
\end{align*}
With our choice of $v(t)=\log|\log|\log t||$,
\begin{align*}
v(r_{i+2})-v(r_i)\sim|\log(i+2)-\log(i)|\sim\frac{2}{i},
\end{align*}
so the sum of $p$-th powers converges to a bound that depends only on $\ell$ and $p$.

\bigskip

3. The third sub-collection consists of balls $B=[a,b]\times\beta$ such that $z_0\notin \ell\beta$. 

Let us study how $r$ varies along $B$. Note that the radius of balls $B$ and $\beta$ equals $\frac{b-a}{2}$. Let $\Delta=\max_\beta d(\cdot,z_0)$ and $\delta=\min_\beta d(\cdot,z_0)$. Then
\begin{align*}
\max_B r=\max\{b,\Delta\},\quad \min_B r=\max\{a,\delta\}.
\end{align*}
Since $z_0\notin\ell\beta$, $\delta\geq (\ell-1)\frac{b-a}{2}$. Then $\Delta\leq\delta+b-a\leq\delta+\frac{2}{\ell-1}\delta=m\delta$. On the other hand, 
\begin{itemize}
  \item either $a<\frac{b}{2}$, $\delta\geq (\ell-1)\frac{b-a}{2}\geq (\ell-1)\frac{b}{4}$,
  \item or $a\geq\frac{b}{2}$.
\end{itemize}
In the first case, 
\begin{align*}
\max\{b,\Delta\}\leq \max\{\frac{4}{\ell-1}\delta,m\delta\}\leq\max\{\frac{4}{\ell-1},m\}\max\{a,\delta\}.
\end{align*}
In the second case,
\begin{align*}
\max\{b,\Delta\}\leq \max\{2a,m\delta\}\leq\max\{2,m\}\max\{a,\delta\}.
\end{align*}
In either case,
\begin{align*}
\max_B r\leq M\min_B r,
\end{align*}
where $M=\max\{m,2,\frac{4}{\ell-1}\}$.

Let $\nu=dt\otimes\mu$ denote $1+Q$-dimensional Hausdorff measure. Set $r_i=M^{-i}$. Let 
\begin{align*}
Y_i=\{(y,z)\in [0,1]\times Z\,;\,r_{i+2}<r(y,z)\leq r_{i}\}, \quad L_i=Lip(v_{|[r_{i+2},r_{i}]}).
\end{align*}
Since $r$ is 1-Lipschitz, $Lip(u_{|Y_i})\leq L_i$.
Each ball $B_j$ of the sub-packing is entirely contained in one of the sets $Y_i$. Therefore the index set is contained in the union of subsets
\begin{align*}
J'_i=\{j\,;\, B_j\subset Y_i\},
\end{align*} 
and for all $p$,
\begin{align*}
\sum_{j}\mathrm{diameter}(w(B_j))^p \leq \sum_{i=0}^{\infty}\sum_{j\in J'_i}\mathrm{diameter}(u(B_j))^p.
\end{align*}
If $j\in J'_i$,
\begin{align*}
\mathrm{diameter}(w(B_{j}))^{1+Q}\leq L_i^{1+Q}\mathrm{diameter}(B_j)^{1+Q}\leq C\,L_i^{1+Q}\nu(B_j).
\end{align*}
Thus
\begin{align*}
\sum_{j\in J'_i}\mathrm{diameter}(u(B_j))^{1+Q}\leq C\,L_i^{1+Q}\nu(Y_i)\leq C'\,(L_i r_{i})^{1+Q},
\end{align*}
since $Y_i\subset B((0,z_0),r_i)$. 

In order to estimate the Lipschitz constant $L_i$, observe that $r$ is 1-Lipschitz. With our choice of $v(t)=\log|\log\|\log t||$, 
\begin{align*}
v'(t)=\frac{1}{t\log t\,\log|\log t|}
\end{align*}
achieves its maximum on $[r_{i+2},r_i]$ at $r_{i+2}=M^{-(i+2)}$, hence
\begin{align*}
L_i\leq\frac{1}{r_{i+2}(i+2)\log(i+2)},\quad L_i r_i\leq\frac{r_i}{r_{i+2}(i+2)\log(i+2)}\leq\frac{M^2}{i+2}.
\end{align*}
This bounds $\sum_{j}\mathrm{diameter}(w(B_j))^{1+Q}$ by a quantity that depends only on $\ell$ and $Q$.

\bigskip

4. At last, we take into account the parameter $\alpha$. When the second factor is the snowflake space $Z^\alpha$, the previous general discussion provides the exponent $1+\frac{Q}{\alpha}$. This can be improved into $\alpha+Q$ for the following reason: on $Z^\alpha$, the distance to $z_0$ is not merely 1-Lipschitz, it is $1/\alpha$-H\"older continuous (although with a constant that deteriorates when getting close to $z_0$). Indeed, if $z,z'\in Z$, with $d(z,z_0)\leq d(z',z_0)$,
\begin{align*}
|d_{Z^\alpha}(z,z_0)-d_{Z^\alpha}(z',z_0)|
&= |d_{Z}(z,z_0)^\alpha-d_{Z}(z',z_0)^\alpha|\\
&\leq \alpha d_{Z}(z,z_0)^{\alpha-1}|d_{Z}(z,z_0)-d_{Z}(z',z_0)|\\
&\leq d_{Z}(z,z_0)^{\alpha-1}d_{Z}(z,z')\\
&= d_{Z^\alpha}(z,z_0)^{\frac{\alpha-1}{\alpha}}d_{Z^\alpha}(z,z')^{1/\alpha}.
\end{align*}

Nothing needs be changed for the first two sub-collections, since an upper bound on energy sums is obtained for any exponent $p>1$. For the third one, consisting of balls $B=[a,b]\times\beta$ such that $z_0\notin\ell\beta$, recall that $\delta\geq(\ell-1)\frac{b-a}{2}$, where
\begin{align*}
\Delta=\max_\beta d_{Z^{\alpha}}(\cdot,z_0)\quad \textrm{and}\quad\delta=\min_\beta d_{Z^{\alpha}}(\cdot,z_0).
\end{align*}
This implies that
\begin{align*}
\delta^{\frac{\alpha-1}{\alpha}}\geq(\ell-1)^{\frac{\alpha-1}{\alpha}}(\frac{b-a}{2})^{\frac{\alpha-1}{\alpha}}
\end{align*}
and
\begin{align*}
\frac{b-a}{2}\leq (\ell-1)^{\frac{1-\alpha}{\alpha}}\delta^{\frac{\alpha-1}{\alpha}}(\frac{b-a}{2})^{1/\alpha}.
\end{align*}
On the other hand, the $1/\alpha$-H\"older character of $d_{Z^{\alpha}}$ leads to
\begin{align*}
\Delta-\delta\leq \delta^{\frac{\alpha-1}{\alpha}}(\frac{b-a}{2})^{1/\alpha}.
\end{align*}
Since
\begin{align*}
\max_B r=\max\{b,\Delta\},\quad \min_B r=\max\{a,\delta\},
\end{align*}
\begin{align*}
\max_B r -\min_B r&\leq \max\{b-a,\Delta-\delta\}\\
&\leq 2(\ell-1)^{\frac{1-\alpha}{\alpha}}\delta^{\frac{\alpha-1}{\alpha}}(\frac{b-a}{2})^{1/\alpha}\\
&\leq M'(\max_B r)^{\frac{\alpha-1}{\alpha}}\mathrm{diameter}(B)^{1/\alpha},
\end{align*}
for $M'=(2/(\ell-1))^{(\alpha-1)/\alpha}$.

Keeping the notation 
\begin{align*}
r_i=M^{-i},\quad Y_i=\{(y,z)\in [0,1]\times Z\,;\,r_{i+2}<r(y,z)\leq r_{i}\},\quad L'_i=Lip(v_{|[r_{i+2},r_{i}]}),
\end{align*}
we see that a ball $B$ contained in $Y_i$ satisfies
\begin{align*}
\mathrm{diameter}(w(B))
&\leq \mathrm{diameter}(v(B))\\
&\leq L'_i(\max_B r -\min_B r)\\
&\leq L'_i M'\,r_{i}^{\frac{\alpha-1}{\alpha}}\mathrm{diameter}(B)^{1/\alpha}.
\end{align*}
Thus
\begin{align*}
\mathrm{diameter}(w(B))^{\alpha+Q}
&\leq {L'_i}^{\alpha+Q} M'^{\alpha+Q}\,r_{i}^{\frac{(\alpha-1)(\alpha+Q)}{\alpha}}\mathrm{diameter}(B)^{1+(Q/\alpha)}\\
&\leq C\,{L'_i}^{\alpha+Q}r_{i}^{\frac{(\alpha-1)(\alpha+Q)}{\alpha}}\nu(B).
\end{align*}
Since $\nu(Y_i)\leq r_i^{1+(Q/\alpha)}$, summing over all balls in the sub-packing contained in $Y_i$ gives
\begin{align*}
\sum_{j\in J'_i}\mathrm{diameter}(w(B))^{\alpha+Q}
&\leq C\,{L'_i}^{\alpha+Q} r_{i}^{\frac{(\alpha-1)(\alpha+Q)}{\alpha}}r_i^{1+(Q/\alpha)}\\
&\leq C\,(L'_i r_i)^{\alpha+Q}.
\end{align*}
The choice of $v(t)=\log|\log|\log t||$ yields again $L_i r_i \leq M^2/(i+2)$, and the sum is bounded above in terms of $\ell$, $\alpha$ and $Q$ only.

5. The final argument, showing that $w\circ\gamma$ has no limit for every proper coarse curve $\gamma$ is the same as in Proposition \ref{parabolicnc}.
\end{pf}

\section{\texorpdfstring{$L^p$}{}-cohomology}
\label{lp}

\subsection{Definition}

Here is one more avatar of the definition of $L^p$ cohomology for metric spaces. This one has the advantage that it does not require any measure. For earlier attempts, see \cite{CT}, \cite{E}.

\begin{defi}
\label{defcoh}
Let $X$ be a metric space. A \emph{$k$-simplex of size $S$} in $X$ is a $k+1$-tuple of points belonging some ball of radius $S$. A \emph{$k$-cochain of size $S$} on $X$ is a real valued function $\kappa$ defined on the set of $k$-simplices of size $S$. Its $L^{p}_{\ell,R,S}$-norm is
\begin{align*}
\n{\kappa}_{L^{p}_{\ell,R,S}}=\sup\{\sum_{j}\sup_{(B_j)^{k+1}}|\kappa|^{p}\,;\,\textrm{all }(\ell,R,S)\textrm{-packings }\{B_j\}\}^{1/p}.
\end{align*}
Let $L^{p}_{\ell,R,S}C^{k}(X)$ denote the space of $k$-cochains with finite $L^{p}_{\ell,R,S}$-norm.

The coboundary operator $d$ maps $k$-cochains to $k+1$-cochains,
\begin{align*}
d\kappa(x_0,\ldots,x_{k+1})&= \kappa(x_1,\ldots,x_{k+1})-\kappa(x_0,x_2,\ldots,x_{k+1})+\cdots\\
&+(-1)^{k+1}\kappa(x_0,\ldots,x_{k}).
\end{align*}
Denote by
\begin{align*}
\mathcal{L}^{q,p}_{\ell,R,S}C^{k}(X)=L^{q}_{\ell,R,S}C^{k}(X)\cap d^{-1}L^p_{\ell,R,S}C^{k+1}(X),
\end{align*}
in order to turn $d$ into a bounded operator $\mathcal{L}^{q,p}_{\ell,R,S}(X)\to L^p_{\ell,R,S}(X)$.
The \emph{$L^{q,p}$-cohomology} of $X$ is 
\begin{align*}
L^{q,p}_{\ell,R,S}H^{k}(X)=\left(\mathrm{ker}(d)\cap L^{p}_{\ell,R,S}C^{k}(X)\right)/d\left(\mathcal{L}^{q,p}_{\ell,R,S}C^{k-1}(X)\right).
\end{align*}
The \emph{exact $L^{q,p}$-cohomology} of $X$ is the kernel of the forgetful map $L^{q,p}_{\ell,R,S}H^{k}(X)\to H^{k}(X)$.

When $p=q$, $\mathcal{L}^{p,p}=L^p$ and $L^{p,p}$-cohomology is simply called $L^p$-cohomology.
\end{defi}

\begin{rem}
The definition of $L^{q,p}_{\ell,0,\infty}H^{k}(X)$ extends to q.s. spaces $X$, and is simply denoted by $L^{q,p}_{\ell}H^{k}(X)$.
\end{rem}

For instance, a function $u:X\to\R$ can be viewed as a $0$-cochain of infinite size, $du(x_1,x_2)=u(x_2)-u(x_1)$ is a $1$-cochain of infinite size, and
\begin{align*}
E^{p}_{\ell,R,S}(u)=\n{du}_{L^{p}_{\ell,R,S}}^{p}.
\end{align*}

\begin{ex}
If $X$ is a compact infinite $d$-Ahlfors regular metric space, then, for all $S>0$, $L^{p}_{\ell,0,S}H^1(X)\not=0$ for $p\geq d$.
\end{ex}
Indeed, in an infinite metric space, one can $\ell$-pack infinitely many small balls. Therefore a function which is $\geq 1$ has infinite $L^{p}_{\ell,0,S}$ norm. Since non constant Lipschitz functions on $X$ have finite energy, and do not belong to any $L^{p}_{\ell,0,S}C^0(X)$, their $L^{p}_{\ell,0,S}$ cohomology classes do not vanish.

\subsection{Link to usual \texorpdfstring{$L^p$}{}-cohomology}

$L^p$-cohomology calculations on manifolds (resp. on simplicial complexes) require the classical de Rham (resp. simplicial) definition of cohomology. There is a de Rham style theorem relating Definition \ref{defcoh} to smooth differential forms (resp. simplicial cochains). It shows up in \cite{genton}. We shall need the more general case of $L^{q,p}$-cohomology, which appears in \cite{ducret}. 

\begin{prop}[\cite{BFP}]
\label{smooth}
Let $1\leq p\leq q<+\infty$. Let $X$ be a bounded geometry Riemannian manifold with boundary or simplicial complex. In the $n$-manifold case, assume that
\begin{align*}
\frac{1}{p}-\frac{1}{q}\leq\frac{1}{n}.
\end{align*}
and $p>1$ if degree $k=n$ (these limitations are absent in the simplicial complex case).
Assume that cohomology of $X$ vanishes uniformly up to degree $k$, i.e. for all $T>0$, there exists $\tilde{T}$ such that for all $x\in X$, the inclusion $B(x,T)\to B(x,\tilde{T})$ induces the 0 map in cohomology up to degree $k$. 

Then for every $R>0$ and $S<+\infty$ and for large enough $\ell\geq 1$, there is a natural isomorphism of $L^{q,p}$-cohomologies $L^{q,p}_{\ell,R,S}H^k(X)\simeq L^{q,p}H^k(X)$. In degree $k+1$, the isomorphism persists provided the space $L^{q,p}H^k(X)$ is replaced with \emph{exact cohomology}, i.e. the kernel $EL^{q,p}H^k(X)$ of the forgetful map $L^{q,p}H^k(X)\to H^k(X)$.
This isomorphism is compatible with multiplicative structures.
\end{prop}

For degree 1 cohomology, the size limit plays no role. Indeed, a 1-cocycle of size $S$ on a simply connected space is the differential of a function, and thus uniquely extends into a 1-cocycle without size limit. The considerations of subsection \ref{radii} show that, \emph{a priori}, all $L^p_{\ell,R,S}$-norms are equivalent on 1-cocycles. For higher degree cohomology, this holds only at the cohomology level, under suitable assumptions, thanks to Proposition \ref{smooth}.

\begin{rem}
\label{lq+ls}
\cite{BFP} contains the following useful remark. On an $n$-manifold, one can define de Rham $L^p$ versus $L^q +L^s$-cohomology. One could do so with cochains, but this would be useless since $\ell^s\subset\ell^q$ when $s\leq q$. Fortunately, under the assumptions of Proposition \ref{smooth}, and if $s\leq q$, this de Rham $L^p$ versus $L^q +L^s$-cohomology is isomorphic to $L^{q,p}_{\ell,R,S}H^\cdot$, and so to $L^{q,p}H^\cdot$.
\end{rem}

\subsection{Functoriality of \texorpdfstring{$L^p$}{}-cohomology}

For every $\ell'\geq 1$, $(R,S,R',S')$-coarse conformal map $f:X\to X'$ induces a bounded linear map 
\begin{align*}
f^*:L^{q,p}_{\ell',R',S'}H^k(X')\to L^{q,p}_{\ell,R,S}H^k(X)
\end{align*}
for suitable $\ell\geq 1$. So, under the assumptions of Proposition \ref{smooth}, a coarse conformal map of a bounded geometry Riemannian manifold with boundary or simplicial complex to a q.s. space $X'$ induces maps
\begin{align*}
f^*:L^{q,p}_{\ell'}H^k(X')\to L^{q,p}H^k(X),
\end{align*}
for all large enough $\ell'$. Exact cohomology is natural as well.

Under the assumptions of Proposition \ref{smooth}, $L^{q,p}$-cohomology is natural under coarse embeddings and a quasi-isometry invariant. On the other hand, it is not clear wether it is natural under large scale conformal maps.

\subsection{Vanishing of 1-cohomology and limits}

\begin{defi}
Let $X$ be a metric space, $Y$ a topological space, and $y\in Y$. Assume $X$ is unbounded. Say a map $f:X\to Y$ tends to $y$ at infinity if for every neighborhood $V$ of $y$, there exists a bounded set $K$ such that $f(x)\in V$ when $x\notin K$.
\end{defi}

\begin{lem}
Let $X$ be an unbounded metric space. Let $q<\infty$. Then every function $u\in L^{q}_{\ell,R,S}C^{0}(X)$ tends to $0$ at infinity. 
\end{lem}

\begin{pf}
Fix $\eps>0$. Let $\{B_j\}$ be an $(\ell,R,S)$-packing such that 
\begin{align*}
\sum_{j}(\sup_{B_j}|u|)^q >\n{u}_{L^{q}_{\ell,R,S}}^q-\eps.
\end{align*}
Pick a finite subfamily which achieves the sum minus $\eps$. The union of this finite subfamily is contained in a ball $K$. If $d(x,\ell K)>\ell R$, add $B(x,R)$ to the finite subfamily to get a larger $(\ell,R,S)$-packing. By definition of energy, $\sup_{B(x,R)}|u|^q<2\eps$. In particular, $|u(x)|<(2\eps)^{1/q}$ outside a bounded set. This shows that $u$ tends to $0$ at infinity.
\end{pf}

\begin{cor}
\label{h=0=>limit}
Let $X$ be an unbounded metric space. Let $q<\infty$. Assume that the $L^{q,p}$-cohomology of $X$ vanishes, i.e. $L^{q,p}_{\ell,R,S}H^{1}(X)=0$. Let $Y$ be a complete metric space. Then every map from $X$ to $Y$ with finite $E^{p}_{\ell,R,S}$ energy has a limit at infinity.
\end{cor}

\begin{pf}
Let $u:X\to Y$ have finite energy. For $y\in Y$, set $v_y(x)=d(u(x),y)$. Then $v_y$ has finite $E^{p}_{\ell,R,S}$ energy. By assumption, there exists a $0$-cochain $w\in L^{q}_{\ell,R,S}C^0(X)$ such that $dw=dv_y$. This implies that $v_y$ has a finite limit $\alpha(y)$ at infinity. $\alpha$ belongs to the closure of the Kuratowski embedding of $Y$ in $L^{\infty}(Y)$. Since $Y$ is complete, the embedding has a closed image, so there exists a point $z\in Y$ such that $\alpha(y)=d(z,y)$ for all $y\in Y$. In particular, $d(u(x),z)=v_z (x)$ tends to $\alpha(z)=0$.\end{pf}

\subsection{Vanishing of reduced 1-cohomology and limits}

\begin{defi}
Let $X$ be a metric space. The \emph{reduced $L^{q,p}$-cohomology} of $X$ is obtained by modding out by the $L^p$-closure of the image of the coboundary operator $d$,
\begin{align*}
L^{q,p}_{\ell,R,S}\bar{H}^{k}(X)=\left(\mathrm{ker}(d)\cap L^{p}_{\ell,R,S}C^{k}(X)\right)/\overline{d\left(\mathcal{L}^{q,p}_{\ell,R,S}C^{k-1}(X)\right)}.
\end{align*}
\end{defi}

\begin{lem}
Let $X$ be an unbounded metric space with a base point. Then for every finite energy function $u$ such that the reduced $L^{q,p}_{\ell,R,S}$-cohomology class of $du$ vanishes, there exists $c\in\R$ such that $u$ converges to $c$ along $p$-almost every based curve. 
\end{lem}

\begin{pf}
Assume that $u_j\in L^{q}_{\ell,R,S}C^0 (X)$ and
\begin{align*}
\n{du_j -du}_{L^{p}_{\ell,R,S}}\textrm{ tends to }0.
\end{align*}
For $t\in\R$, let $\Gamma_{t,+}$ (resp. $\Gamma_{t,-}$) be the family of based curves $\gamma$ along which $u$ has a finite limit and $\lim u\circ\gamma\geq t$ (resp. $\leq t$). Fix $s<t$. Let $v_j =\frac{2}{t-s}(u-u_j)$. 

Assume that there exists a based curve $\gamma\in\Gamma_{t,+}$ such that, for infinitely many $j$, $\mathrm{length}(v_j\circ\gamma)\leq 1$. For those $j$'s, for every $\gamma'\in\Gamma_{s,-}$, $\mathrm{length}(v_j\circ\gamma)\leq 1$. Indeed, along the bi-infinite curve obtained by concatenating $\gamma$ and $\gamma'$, the total variation of $v_j$ is $\geq 2$. Therefore, for infinitely many $j$'s,
\begin{align*}
\mathrm{mod}_{p,\ell,R,S}(\Gamma_{s,-})\leq E(v_j)=(\frac{2}{t-s} \n{du_j -du}_{L^{p}_{\ell,R,S}})^{p},
\end{align*}
and $\mathrm{mod}_{p,\ell,R,S}(\Gamma_{s,-})=0$. 

Otherwise, for each $\gamma\in\Gamma_{t,+}$, for all but finitely many $j$'s, \break $\mathrm{length}(v_j\circ\gamma)\geq 1$. $\Gamma_{t,+}$ is the union of subfamilies
\begin{align*}
\Gamma_{t,+,J}=\{\gamma\in\Gamma_{t,+}\,;\,\forall j\geq J,\,\mathrm{length}(v_j\circ\gamma)\geq 1\},
\end{align*}
each of which has vanishing modulus. By stability under countable unions, $\mathrm{mod}_{p,\ell,R,S}(\Gamma_{t,+})=0$. 

Let $c$ be the supremum of all $t\in\R$ such that $\mathrm{mod}_{p,\ell,R,S}(\Gamma_{t,+})>0$. By stability under countable unions, the family of based curves along which $u$ has a finite limit $>c$ has vanishing modulus. If $c=-\infty$, for every $n\in\Z$, the family of based curves along which $u$ has a finite limit $\leq n$ has vanishing modulus. Thus the family of all based curves has vanishing modulus, and the Lemma is proved. Otherwise, $\mathrm{mod}_{p,\ell,R,S}(\Gamma_{s,-})=0$ for all $s<c$. By stability under countable unions, the family of based curves along which $u$ has a finite limit $<c$ has vanishing modulus. This shows that
$u$ tends to $c$ along almost every based curve. \end{pf}

\begin{cor}
\label{cohomvanishing}
Let $X$ be an unbounded metric space. Let $Y$ be a complete metric space. Let $u:X\to Y$ have finite $E^{p}_{\ell,R,S}$ energy. Assume that the reduced $L^{q,p}$-cohomology of $X$ vanishes, i.e. $L^{q,p}_{\ell,R,S}\bar{H}^{1}(X)=0$. Then $u$ has a common limit along $p$-almost every based curve.
\end{cor}

\begin{pf}
For $y\in Y$, set $v_y(x)=d(u(x),y)$. Then $v_y$ has finite $E^{p}_{\ell,R,S}$ energy. By assumption, $dv_y$ belongs to the $L^{p}_{\ell,R,S}$-closure of $d\mathcal{L}^{q,p}_{\ell,R,S}C^0(X)$. This implies that $v_y$ has a finite limit $\alpha(y)$ along almost every based curve. $\alpha$ belongs to the closure of the Kuratowski embedding of $Y$ in $L^{\infty}(Y)$. Since $Y$ is complete, the embedding has a closed image, so there exists a point $z\in Y$ such that $\alpha(y)=d(z,y)$ for all $y\in Y$. In particular, $d(u(x),z)=v_z (x)$ tends to $\alpha(z)=0$ along almost every based curve.\end{pf}

\subsection{\texorpdfstring{$p$}{}-separability}

\begin{defi}
Let $X$ be a q.s. space. Let $\mathcal{E}_{p,\ell}$ denote the space of continuous functions on $X$ with finite $(p,\ell)$-energy.
\end{defi}

\begin{defi}
Say a q.s. space $X$ is $p$-separated if for every large enough $\ell$, $\mathcal{E}_{p,\ell}$ separates points and complements of points are $(p,\ell)$-parabolic. If $X$ is non-compact, one requires further that $X$ itself be $(p,\ell)$-parabolic.
\end{defi}

\begin{ex}
\label{expsep}
$Q$-Ahlfors-regular metric spaces are $p$-separated for all $p\geq Q$.
\end{ex}

\begin{pf}
Proposition \ref{lipenergy} shows that Lipschitz functions with bounded support have finite energy. They separate points. Propositions \ref{parabolicnc} and \ref{parabolic} establish parabolicity of $X$ and of point complements.\end{pf}

\begin{prop}
\label{separated}
Let $1\leq p,q<+\infty$. Let $X$ be a locally compact metric space. Let $X'$ be a locally compact $p$-separated q.s. space. Let $f:X\to X'$ be a coarse conformal map. Then there exists $R>0$ such that for all $S\geq R$ and $\ell'>1$, for all large enough $\ell$, \begin{itemize}
  \item either the induced map
\begin{align*}
f^*:EL^{q,p}_{\ell'}\bar{H}^{1}(X')\to EL^{q,p}_{\ell,R,S}\bar{H}^{1}(X)
\end{align*}
in reduced exact $L^{q,p}$-cohomology does not vanish, 
  \item or $X$ is $(p,\ell,R,S)$-parabolic.
\end{itemize}
\end{prop}

\begin{pf}
Let us treat first the simpler case when unreduced $L^p$-cohomo\-logy vanishes. Let $f:X\to X'$ be a coarse conformal map. Assume that $f$ has distinct accumulation points $x'_1$ and $x'_2$ at infinity. By assumption, there exists a continuous function $v$ on $X'$ with finite $(p,\ell')$-energy such that $v(x'_1)\not=v(x'_2)$. Then, for all $R\leq S$, $v\circ f$ has finite $E^{p}_{\ell,R,S}$ energy for suitable $\ell$. If $L^{q,p}_{\ell,R,S}H^{1}(X)$ vanishes, according to Corollary \ref{h=0=>limit}, $v\circ f$ has a limit at infinity. This should be at the same time $v(x'_1)$ and $v(x'_2)$, contradiction. We conclude that $f$ has at most one accumulation point at infinity. Hence either it has a limit $x'$, or it tends to infinity.

In either case, there exists a finite energy function $w:X'\to\R$ that has no limit along $p$-almost every coarse curve converging to $x'$ (resp. to infinity). Then $w\circ f$ has finite $E^{p}_{\ell,R,S}$ energy as well, it must have a finite limit in $\R$. This contradicts the fact that the family of based $(1,1)$-curves in $X$ has positive $(p,\ell,R,S)$-modulus. We conclude that $f^*$ does not vanish on $L^{q,p}_{\ell}H^{1}(X')$.

Assume that $f$ induces a trivial map in reduced cohomology and that $X$ is non-$(p,\ell,R,S)$-parabolic. Equip $\mathcal{E}_{p,\ell}$ with the topology of uniform convergence on compact sets. Let $D\subset \mathcal{E}_{p,\ell}$ be a countable dense subset. Then $D$ still separates points. We know that for all $v\in D$, $v\circ f\circ \gamma$ has a common limit $y_v$ for almost every based $(1,1)$-curve $\gamma$ in $X$. For $v\in D$, let $\Gamma_v$ be the family of $(R,S)$-based curves $\gamma\subset X$ such that $v\circ f\circ \gamma$ does not have a limit or has a limit which differs from $y_v$. Then $\Gamma=\bigcup_{v\in D}\Gamma_v$ has vanishing $(p,\ell,R,S)$-modulus. Let $\Gamma'$ be the complementary family. Since $X$ is non-$(p,\ell,R,S)$-parabolic, $\Gamma'$ is non-empty. Fix two based curves $\gamma,\gamma'\in\Gamma'$. Assume that $f\circ\gamma$ and $f\circ\gamma'$ have distinct accumulation points $x'_1$ and $x'_2$ in $X'$. Let $v\in D$ be such that $v(x'_1)\not= v(x'_2)$. By construction, $v\circ\gamma$ and $v\circ\gamma'$ converge to $y_v$. Since $v$ is continuous, $v\circ\gamma$ subconverges to $v(x'_1)$ and $v\circ\gamma'$ to $v(x'_2)$, contradiction. We conclude that $f$ has a common limit $x'$ along all $\gamma\in\Gamma'$. The argument ends in the same manner.\end{pf}

\subsection{Relative \texorpdfstring{$p$}{}-separability}

Here comes a relative version of Proposition \ref{separated}, motivated by the case of warped products.

\begin{defi}
Let $X$ be a q.s. space. Let $u:X\to Y$ be a continuous map to a topological space. Say $X$ is \emph{$p$-separated relatively to $u$} if for every large enough $\ell$, 
\begin{enumerate}
  \item $\mathcal{E}_{p,\ell}\cup\{u\}$ separates points.
  \item Complements of points in $X$ are $(p,\ell)$-parabolic.
\end{enumerate}
\end{defi}

\begin{ex}
\label{exrelsep}
Let $Z$ be a compact $Q$-Ahlfors-regular metric space. Let $0<\alpha\leq 1$. Let $Z^\alpha=(Z,d_Z^\alpha)$ be a snowflaked copy of $Z$. Let $X=\mathbb{D}\times Z^\alpha$. Let $u$ be the projection to the first factor $Y=[0,1]$. Then $X$ is $p$-separated relative to $u$ for all $p\geq \alpha+Q$.
\end{ex}
Indeed, the projection to the second factor equipped with $d_Z$ has finite $(p,\ell)$-energy for all $p\geq\alpha+Q$ and all $\ell>1$, see Example \ref{extwisted}. Together with $u$, it separates points. Proposition \ref{twisted} states that complements of points in $X$ are $(p,\ell)$-parabolic. 

\begin{prop}
\label{relseparated}
Let $1\leq p,q<+\infty$. Let $X$ be a locally compact metric space containing at least one based $(1,1)$-curve. Let $X'$ be a compact q.s. space which is $p$-separated relatively to a map $u:X'\to Y$. Let $f:X\to X'$ be a coarse conformal map. Assume that the map $u\circ f$ tends to some point $y\in Y$ at infinity.
Then there exists $R>0$ such that for all $S\geq R$ and $\ell'>1$, for all large enough $\ell$,
\begin{itemize}
  \item either the induced map
\begin{align*}
f^*:EL^{q,p}_{\ell'}\bar{H}^{1}(X')\to EL^{q,p}_{\ell,R,S}\bar{H}^{1}(X)
\end{align*}
in reduced exact $L^{q,p}$-cohomology does not vanish, 
  \item or $X$ is $(p,\ell,R,S)$-parabolic.
\end{itemize}
\end{prop}

\begin{pf}
By contradiction. Assume that $f$ induces a trivial map in reduced cohomology and that $X$ is non-$(p,\ell,R,S)$-parabolic. Equip $\mathcal{E}_{p,\ell}$ with the topology of uniform convergence. Let $D\subset \mathcal{E}_{p,\ell}$ be a countable dense subset. Then $D\cup\{u\}$ still separates points. For all $v\in D$, $v\circ f\circ \gamma$ has a common limit $t_v$ for almost every based $(1,1)$-curve $\gamma$ in $X$. For $v\in D$, let $\Gamma_v$ be the family of $(R,S)$-based curves $\gamma\subset X$ such that $v\circ f\circ \gamma$ does not have a limit or has a limit which differs from $t_v$. Then $\Gamma=\bigcup_{v\in D}\Gamma_v$ has vanishing $(p,\ell,R,S)$-modulus. Let $\Gamma'$ be the complementary family. Since $X$ is non-$(p,\ell,R,S)$-parabolic, $\Gamma'$ is non-empty. Fix two based curves $\gamma,\gamma'\in\Gamma'$. Assume that $f\circ\gamma$ and $f\circ\gamma'$ have distinct accumulation points $x'_1$ and $x'_2$ in $X'$. Since $u(x'_1)=u(x'_2)=y$, there exists $v\in D$ such that $v(x'_1)\not= v(x'_2)$. By construction, $v\circ\gamma$ and $v\circ\gamma'$ converge to $t_v$. Since $v$ is continuous, $v\circ\gamma$ subconverges to $v(x'_1)$ and $v\circ\gamma'$ to $v(x'_2)$, contradiction. We conclude that $f$ has a common limit $x'$ along all $\gamma\in\Gamma'$. 

Let $w:X'\to\R$ be a finite energy function that has no limit along $(p,\ell)$-almost every coarse curve converging to $x'$. Then $w\circ f$ has finite $E^{p}_{\ell,R,S}$ energy as well, it must have a finite limit in $\R$. This contradicts the fact that the family of based $(1,1)$-curves in $X$ has positive $(p,\ell,R,S)$-modulus. We conclude that either $f^*$ does not vanish on $L^{q,p}_{\ell}H^{1}(X')$ or $X$ is $p$-parabolic.
\end{pf}

\section{Lack of coarse conformal maps}
\label{lack}

\begin{thm}
\label{I}
Let $1<p,q<+\infty$. Let $X$ be a Riemannian manifold or a simplicial complex with bounded geometry. %and uniform vanishing of 1-cohomology. 
In the $n$-manifold case, assume further that $0\leq\frac{1}{p}-\frac{1}{q}\leq\frac{1}{n}$. Assume that $X$ is non-$p$-parabolic and that $EL^{q,p}\bar{H}^{1}(X)=0$. 
\begin{enumerate}
  \item Let $X'$ be a $p$-Ahlfors-regular metric space. Then there can be no coarse conformal maps $X\to X'$. 
  \item Let $0<\alpha\leq 1$. Let $X'$ be a warped product $\mathbb{D}\times Z^{\alpha}$ where $Z$ is a compact $p-\alpha$-Ahlfors-regular metric space. For every coarse conformal map $X\to X'$,  the projected map to the first factor $X\to [0,1]$ cannot tend to 0. 
\end{enumerate}
\end{thm}

\begin{pf} It follows from Examples \ref{expsep} and \ref{exrelsep}, Propositions \ref{separated}, \ref{relseparated} and \ref{smooth}. \end{pf}

\subsection{Examples}

For nilpotent groups, reduced $L^p$-cohomology vanishes. Indeed, such groups admit unbounded central subgroups. A central element in $G$ acts by a translation of $G$, i.e. moves points a bounded distance away. The Corollary on page 221 of \cite{G} applies: reduced $L^p$-cohomology vanishes in all degrees, in particular in degree 1. 
 
A nilpotent group of homogeneous dimension $Q$ is $p$-parabolic if and only if $p\geq Q$. $(p,\ell,R,\infty)$-parabolicity for $p\geq Q$, $R>0$ and $\ell>1$ follows from the asymptotics of volume of balls, \cite{growth}, and Remark \ref{lsparabolic}.
Non-$(p,\ell,R,S)$-parabolicity for $p<Q$ will be proved below, in Proposition \ref{isodimgps} and Corollary \ref{isop=>snp}. Carnot groups in their Carnot-Carath\'eodory metrics are even $(p,\ell)$-parabolic for $p\geq Q$ and $\ell>1$, according to Proposition \ref{parabolicnc}, since they are $Q$-Ahlfors regular.

\medskip

Non-elementary hyperbolic groups, \cite{GH}, have infinite isoperimetric dimension, hence they are never $p$-parabolic (again, this follows from Proposition \ref{isodimgps} and Corollary \ref{isop=>snp}). Their $L^p$-cohomology vanishes for $p$ in an interval starting from 1, whose upper bound is denoted by $\mathrm{CohDim}$ (\cite{BP}).
 
Conformal dimension $\mathrm{ConfDim}$ arises as the infimal Hausdorff dimension of Ahlfors-regular metrics in the quasi-symmetric gauge of the ideal boundary. By definition, the quasi-symmetric gauge is the set of metrics which are quasi-symmetric to a visual quasi-metric (all such quasi-metrics are mutually quasi-symmetric), \cite{BP}. Quite a number of results on $\mathrm{CohDim}$ and $\mathrm{ConfDim}$ can be found in recent works by Marc Bourdon, \cite{BnK}, John Mackay \cite{McK} and their co-authors.

\subsection{Maps to nilpotent groups}
\label{tonil}

\begin{cor}
\label{carnot}
Let $G$ and $G'$ be nilpotent Lie group of homogeneous dimensions $Q$ and $Q'$. Assume that $G'$ is Carnot and equipped with a homogeneous Carnot-Carath\'eodory metric. If there exists a coarse conformal map $G\to G'$, then $Q\leq Q'$.
\end{cor}

\begin{pf}
$G'$, a Carnot group in its Carnot-Carath\'eodory metric, is $Q'$-Ahlfors-regular. Since reduced $L^p$-cohomology vanishes, Theorem \ref{I} forbids existence of a coarse conformal map $G\to G'$ unless $G$ is $Q'$-parabolic. This implies that $Q\leq Q'$. 
\end{pf}
Note that no properness assumption was made. Also, a stronger result will be obtained by a different method in Corollary \ref{arbtonil}.

\begin{cor}
\label{niltonil}
Let $G$ be a finitely generated group. Let $G'$ be a nilpotent Lie or finitely generated group. If there exists a uniformly conformal map $G\to G'$, then $G$ is itself virtually nilpotent, and $d(G)\leq d(G')$.
\end{cor}

\begin{pf}
It is $Q'$-Ahlfors regularity of $G'$ \textrm{in the large} which is used here, and Corollaries \ref{lsc=>para} and \ref{parabolicinthe large}. Indeed, $G'$ is $p$-parabolic for $p=d(G')$, therefore so is $G$. Proposition \ref{lp=>snp} implies that the isoperimetric dimension of $G'$ is at most $p$. Proposition \ref{isodimgps} tells us that $G$ must be virtually nilpotent and $d(G)\leq p=d(G')$.
\end{pf}
 
\begin{cor}
\label{hyptocar}
Let $G$ be a non-elementary hyperbolic group. Let $G'$ be a Carnot group of homogeneous dimension $Q'$ equipped with its Carnot-Carath\'e\-odory metric. If there exists a coarse conformal map $G\to G'$, then $\mathrm{CohDim}(G)\leq Q'$.
\end{cor}

\begin{pf}
Since $G'$ is $Q'$-Ahlfors regular and non-elementary hyperbolic groups are never $p$-parabolic, Theorem \ref{I} provides this upper bound on $\mathrm{CohDim}(G)$.
\end{pf}

\begin{cor}
\label{hyptonil}
Let $G$ be a non-elementary hyperbolic group. Let $G'$ be a nilpotent group of homogeneous dimension $Q'$. If there exists a uniformly conformal map $G\to G'$, then $\mathrm{CohDim}(G)\leq Q'$.
\end{cor}

\begin{pf}
Corollaries \ref{lsc=>para} and \ref{parabolicinthe large} apply, since non-elementary hyperbolic groups are never $p$-parabolic.
\end{pf}

\subsection{Maps to hyperbolic groups}
\label{tohyp}

\begin{cor}
\label{hyptohyp}
Let $G$, $G'$ be non-elementary hyperbolic groups. If there exists a uniformly conformal map $G\to G'$, then 
$$\mathrm{CohDim}(G)\leq \mathrm{ConfDim}(G').$$
\end{cor}

\begin{pf}
Let $d'$ be an Ahlfors-regular metric in the gauge of $\partial G'$, of Hausdorff dimension $Q'$. According to \cite{carrasco}, there exist a bounded geometry hyperbolic graph $X$, a visual quasi-metric $d_o$ on $\partial X$ and a bi-Lipschitz homeomorphism $q:(\partial G',d')\to(\partial X,d_o)$, arising from a quasi-isometry $q:G'\to X$. Set $Z=(\partial X,(q^{-1})^*d')$. Pick $0<\alpha\leq 1$. Let $X'=\mathbb{D}\times Z^{\alpha}$. According to Proposition \ref{poinc=>rc}, the Poincar\'e model of $X$ is a roughly conformal map $\pi:X\to X'$. If $f:G\to G'$ is uniformly conformal, then $f'=\pi\circ q\circ f:G\to X'$ is coarsely conformal (Proposition \ref{composition}). Furthermore, $q\circ f$ is proper, so the projection of $f'$ to the first factor tends to $0$. 

Since $G$ is never $p$-parabolic, Theorem \ref{I} asserts that $\mathrm{CohDim}(G)\leq\alpha+Q'$. Taking the infimum over $\alpha\in(0,1)$ and $Q'\geq \mathrm{ConfDim}(G)$, we get $\mathrm{CohDim}(G)\leq \mathrm{ConfDim}(G')$. 
\end{pf}

\begin{ex}
Fuchsian buildings.
\end{ex}
Right-angled Fuchsian buildings (also known as Bourdon buildings) $X_{p,q}$ are universal covers of orbihedra having one $p$-sided polygon, $p$ even, with trivial face group, cyclic $\Z/q\Z$ edge groups and direct product $\Z/q\Z\times \Z/q\Z$ vertex groups. The conformal dimension and the cohomological dimension of $X_{p,q}$ are both equal to $1+\frac{\log q}{\arg\cosh(\frac{p-2}{2})}$, \cite{Bo}, \cite{BP}. As $p$ and $q$ vary, these numbers fill a dense subset of $[1,+\infty)$.

There are obvious isometric embeddings $X_{p,q}\to X_{p,q'>q}$ and a jungle of bi-Lipschitz embeddings $X_{2p-4,q}\to X_{p,q}$, $X_{3p-8,q}\to X_{p,q}$,.... The only known restriction on the existence of uniform/coarse embeddings $X_{p,q}\to X_{p',q'}$ is provided by Corollary \ref{hyptohyp} or, alternatively, by D. Hume, J. Mackay and R. Tessera's $p$-separation estimates, \cite{HMT}.

\begin{cor}
\label{niltohyp}
Let $G$ be a nilpotent Lie group of homogeneous dimension $Q$. Let $G'$ be a hyperbolic group. If there exists a uniformly conformal map $G\to G'$, then $Q\leq \mathrm{ConfDim}(G')$.
\end{cor}

\begin{pf}
Only the last paragraph of the proof of Corollary \ref{hyptohyp} needs be changed. In this case, reduced $L^p$-cohomology vanishes always, Theorem \ref{I} asserts that $G$ must be $p$-parabolic for $p=\alpha+Q'$, hence $Q\leq \alpha+Q'$. Taking the infimum over $\alpha\in(0,1)$ and $Q'$, we get $Q\leq \mathrm{ConfDim}(G')$. \end{pf}

\begin{ex}
This is sharp. For instance, the uniform embeddings $\R^{n-1}\to H^n_\R$ and $Heis^{2m-1}\to H^m_\mathbb{R}$ are uniformly conformal.
\end{ex}
More generally, every Carnot group $G$ is a subgroup of the hyperbolic Lie group $G'=\R\ltimes G$, where $\R$ acts on $G$ through Carnot dilations. The homogeneous dimension of $G$ is equal to the conformal dimension of $G'$, \cite{Pa}. This provides a uniformly conformal map $G\to G'$, according to Lemma \ref{coarse=>lsc}.

\subsection{Proof of Theorem \ref{dd} and Corollary \ref{udd}}

Theorem \ref{dd} is a combination of Corollaries \ref{niltonil}, \ref{hyptonil}, \ref{hyptohyp} and \ref{niltohyp} applied to the subclass of large scale conformal maps. Corollary \ref{udd} is the special case of uniform/coarse embeddings.

\section{Large scale conformal isomorphisms}
\label{iso}

\subsection{Capacities}
\label{capacity}

\begin{defi}
Let $X$ be a metric space, let $K\subset X$ be a bounded set. The \emph{$(p,\ell,R,S)$-capacity} of $K$, $\mathrm{cap}_{p,\ell,R,S}(K)$, is the infimum of $E^{p}_{\ell,R,S}$-energies of functions $u:X\to [0,1]$ which take value 1 on $K$ and have bounded support.
\end{defi}

\begin{rem}
If $\mathrm{cap}_{p,\ell,R,S}(\{o\})=0$, then $X$ is $(p,\ell,R,S)$-parabolic.
\end{rem}

\begin{pf}
The capacity of the one point set $\{o\}$ bounds from above the  $(p,\ell$, $R,S)$-modulus of the family of all $(1,1)$-curves based at $o$.
\end{pf}

Note that $\mathrm{cap}_{p,\ell,R,S}(\{o\})=0$ implies that for every bounded set $K$, $\mathrm{cap}_{p,\ell',R',S'}(K)=0$ for suitable constants. 

\begin{prop}
Let $X$ and $X'$ be locally compact, noncompact metric spaces. Let $f:X\to X'$ be a large scale conformal map. Then, for every $R'>0$, there exists $R>0$ and for every $\ell'\geq 1$, there exist $\ell\geq 1$ and $N'$ such that, for all compact sets $K\subset X$,
\begin{align*}
\mathrm{cap}_{p,\ell,R,\infty}(K)\leq N'\,\mathrm{cap}_{p,\ell',R',\infty}(f(K)).
\end{align*}
\end{prop}

\begin{pf}
If $u:X'\to[0,1]$ has compact support and $u(f(K))=1$, then $u\circ f$ has compact support and $(u\circ f)(K)\geq 1$, thus $E^{p}_{\ell,R,\infty}(u\circ f)\geq\mathrm{cap}_{p,\ell,R,\infty}(K)$. We know from Lemma \ref{energytransport} that
\begin{align*}
E^{p}_{\ell,R,\infty}(u\circ f)\leq N'\,E^{p}_{\ell',R',\infty}(u).
\end{align*}
Taking the infimum over all such functions $u$,
\begin{align*}
\mathrm{cap}_{p,\ell,R,\infty}(K)\leq N'\,\mathrm{cap}_{p,\ell',R',\infty}(f(K)).
\end{align*}
\end{pf}

\subsection{Non-parabolicity and \texorpdfstring{$L^{q,p}$}{} cohomology}

\begin{defi}
A metric space $X$ is \emph{uniformly perfect in the large} if there exists a constant $c>0$ such that, for all $x\in X$ and large enough $T$, $B(x,T)\setminus B(x,cT)\not=\emptyset$.
\end{defi}
Unbounded geodesic spaces (e.g. graphs, Riemannian manifolds), and spaces roughly isometric to such (e.g. locally compact groups) are uniformly perfect in the large. The point of this property is to ensure that $R$-volumes (meaning the number of disjoint $R$-balls that one can pack inside) of large balls are large.

\begin{lem}
Let $X$ be a metric space which is uniformly perfect in the large. Fix a radius $R>0$. Let $\mathrm{vol}_R(B)$ denote the maximal number of disjoint $R$-balls that can be packed in $B$, and $v_R(T)=\inf_{x\in X} \mathrm{vol}_R(B(x,T))$. Then $v_R(T)$ tends to infinity with $T$.
\end{lem}

\begin{pf}
Given a ball $B=B(x_0,T)$, uniform perfectness, applied in \break$B(x_0,\frac{2T}{2+c})$, provides a point $x_1\in B(x_0,\frac{2T}{2+c})\setminus B(x,\frac{2cT}{2+c})$. Then \break $B(x_1,\frac{c}{2+c}T)\subset B(x,T)\setminus B(x_0,\frac{c}{2+c}T)$. Iterating the construction produces a sequence of disjoint balls $B(x_j,(\frac{c}{2+c})^{j+1} T)$ in $B(x,T)$.  If $n=\lfloor\log_{c/2+c}(T/R)\rfloor$, we get $n$ disjoint $R$-balls in $B(x,T)$.\end{pf}

\begin{lem}
\label{lpco=>snp}
Let $X$ be a metric space which is uniformly perfect in the large. Fix $\ell\geq 1$ and $S\geq R>0$. If $L^{q,p}_{\ell,R,S}H^1(X)=L^{q,p}_{\ell,R,S}\bar{H}^1(X)$ for some finite $q$, then the capacity of balls tends to infinity uniformly with their radius: there exists a function $\kappa_{\ell,R}$ such that $\kappa_{\ell,R}(T)$ tends to infinity as $T\to\infty$, and such that for every ball $B(x,T)$ of radius $T$,
\begin{align*}
\mathrm{cap}^{p}_{\ell,R,S}(B(x,T))\geq\kappa_{\ell,R}(T).
\end{align*}

In particular, $X$ is non-$(p,\ell,R,S)$-parabolic. 
\end{lem}

\begin{pf}
By assumption, the coboundary $d:\mathcal{L}^{q,p}_{\ell,R,S}C^{0}(X)\to L^{p}_{\ell,R,S}C^{1}(X)$ has a closed image. Its kernel consists of constant functions, which can be modded out. $d$ becomes a continuous isomorphism between Banach spaces. According to the isomorphism theorem, $d$ has a bounded inverse. Thus there exists a constant $C$ such that, for every function $u\in L^{q}_{\ell,R,S}C^{0}(X)$, there exists a constant $c_u$ such that 
\begin{align*}
\n{u-c_u}_{L^{p}_{\ell,R,S}C^0(X)}\leq C\,\n{du}_{L^{p}_{\ell,R,S}C^{1}(X)}.
\end{align*}
Since $vol_R(X)$ is infinite, constants do not belong to $L^q_{\ell,R,S}C^0(X)$, so $c_u=0$. For functions $u:X\to [0,1]$ with bounded support, this translates into
\begin{align*}
\n{u}_{L^{q}_{\ell,R,S}C^0(X)}\leq C\,E^{p}_{\ell,R,S}(u)^{1/p}.
\end{align*}
Let $B(x,T)$ be a large ball. Uniform perfectness ensures that a logarithmic number of disjoint balls of radius $\ell R$ can be packed into $B(x,T)$. If $u=1$ on $B(x,T)$, using this packing, we get a lower bound on $\n{u}_{L^{q}_{\ell,R,S}C^0}^q$ of the order of $\log(T/\ell R)$ which depends only on $R$, $\ell$ and $T$. This shows that the $p$-capacity of $B(x,T)$ tends to infinity with $T$.\end{pf}

\subsection{\texorpdfstring{$L^{q,p}$}{}-cohomology and isoperimetric dimension}\label{isop}

\begin{defi}
Let $X$ be a Riemannian manifold. Say that $X$ has isoperimetric dimension $\geq d$ if compact subsets $D\subset X$ with smooth boundary and sufficiently large volume satisfy
\begin{align*}
\mathrm{volume}(D)\leq C\,\mathrm{volume}(\partial D)^{\frac{d}{d-1}}.
\end{align*}
\end{defi}

\begin{lem}
\label{isopdim=>redcoh}
Let $X$ be a Riemannian $n$-manifold with bounded geometry. If $X$ has isoperimetric dimension $\geq d>1$, then, for $p<d$, $EL^{q,p}H^1(X)=EL^{q,p}\bar{H}^1(X)$ for all $1\leq p< d$ and $q<\infty$ such that $\frac{1}{p}-\frac{1}{q}=\frac{1}{\max\{n,d\}}$.
\end{lem}

\begin{pf}
The isoperimetric profile of a Riemannian manifold $X$ is
\begin{align*}
I_X(v)=\inf\{\mathrm{volume}(\partial D)\,;\,\mathrm{volume}(D)=v\}.
\end{align*}
For instance, if $X=\R^n$ is Euclidean $n$-space,
\begin{align*}
I_{\R^n}(v)=c_n v^{\frac{n-1}{n}}.
\end{align*}
According to \cite{Munoz-Flores Nardulli}, the isoperimetric profile of a bounded geometry Riemannian manifold is continuous.
%http://cvgmt.sns.it/paper/2408/
According to \cite{Berard-Meyer}, for a bounded geometry Riemannian $n$-manifold, for every $C<c_n$, there exists $v_0$ such that 
\begin{align*}
I_X(v)\geq C\,v^{\frac{n-1}{n}}\quad \textrm{if}\quad v\leq v_0.
\end{align*}
By assumption, there exists a $c>0$ such that $I_X(v)\geq c\,v^{\frac{d-1}{d}}$ for large $v$. Up to reducing the constant $c$, we can assume that this estimate holds as soon as $v\geq v_0$.

Following a classical argument, let us check that an $L^1$ Sobolev inequality holds. Let $u$ be a smooth compactly supported function on $X$. Let $u_t$ be the indicator function of the superlevel set $\{|u|>t\}$, i.e. $u_t(x)=1$ if $|u(x)|>t$, $u_t(x)=0$ otherwise. Let $t_0=\sup\{t\geq 0\,;\,\mathrm{volume}(\{|u|>t\})\geq v_0$. Let $u'$ be the function that is equal to $|u|$ on $\{|u|\leq t_0\}$, and to $t_0$ elsewhere. Let $u''=|u|-u'$. Then
\begin{align*}
\|d|u|\|_1=\|du'\|_1+\|du''\|_1.
\end{align*}
and
\begin{align*}
u'=\int_{0}^{t_0}u_t\,dt,\quad u''=\int_{t_0}^{+\infty}u_t\,dt.
\end{align*}
Using the isoperimetric inequality for volumes $\geq v_0$, we estimate the $L^{d'}$-norm of $u'$, where $d'=\frac{d}{d-1}$,
\begin{align*}
\|u'\|_{d'}&\leq \int_{0}^{t_0}\|u_t\|_{d'}\,dt\\
&= \int_{0}^{t_0}(\mathrm{volume}(\{|u|>t\})^{1/d'}\,dt\\
&\leq \frac{1}{c}\int_{0}^{t_0}(\mathrm{volume}(\{|u|=t\})\,dt\\
&= \frac{1}{c}\|du'\|_1.
\end{align*}
A similar estimate, using the isoperimetric inequality for volumes $\leq v_0$, gives
\begin{align*}
\|u''\|_{n'}&\leq&\frac{1}{C} \|du''\|_1.
\end{align*}
If $d\geq n$, since $u''$ vanishes outside a set of volume $\leq v_0$, H\"older's inequality gives
\begin{align*}
\|u''\|_{d'}\leq v_0^{\frac{1}{n}-\frac{1}{d}}\|u''\|_{n'},
\end{align*}
hence 
\begin{align}\label{sobolev}
\|u\|_{d'}\leq C\,\|du\|_1.
\end{align} 
It follows that $\|u\|_{q}$ is controlled by $\|du\|_p$ provided $p<d$ and $\frac{1}{p}-\frac{1}{q}=\frac{1}{d}$. Indeed, replacing $|u|$ with $|u|^r$, $r\geq 1$ in inequality (\ref{sobolev}), and applying H\"older's inequality, one finds that
\begin{align*}
(\int |u|^{rd'})^{1/d'} \leq C\,\int |u|^{r-1}|du'|\leq C\,(\int |u|^{rd'})^{\frac{r-1}{rd'}}(\int |du|^{\frac{rd'}{rd'-r+1}})^{\frac{rd'-r+1}{rd'}}.
\end{align*}
If $p<d$, on can pick $q$ such that $\frac{1}{p}-\frac{1}{q}=\frac{1}{d}$ and $r=q/d'\geq 1$. Then
\begin{align*}
\|u\|_{q}\leq C\,\|du\|_{p}.
\end{align*}

If $d<n$, one more step is needed. Note that, for any $r\geq 1$, the decomposition of $v=|u|^r$ is $v'=|u'|^r$ and $v''=|u''|^r$. Therefore one can apply the above method to $|u'|^{r}$ and $|u''|^r$ respectively.  For $u'$, one can use $q$ such that $\frac{1}{p}-\frac{1}{q}=\frac{1}{d}$ and get
\begin{align*}
\|u'\|_{q}\leq C\,\|du'\|_{p}.
\end{align*}
For $u''$, since $p<d<n$, one can pick $s$ such that $\frac{1}{p}-\frac{1}{s}=\frac{1}{n}$ and get
\begin{align*}
\|u''\|_{s}\leq C\,\|du''\|_{p}.
\end{align*}
Summing up, since $\|du\|_{p}^p=\|du'\|_{p}^p+\|du''\|_{p}^p$, one finds that
\begin{align*}
\||u|\|_{L^{s}+L^{q}}\leq\|du\|_p.
\end{align*}
This says that reduced and unreduced de Rham $L^p$ versus $L^q +L^s$-cohomology coincide. According to Remark \ref{lq+ls}, this is equivalent to
\begin{align*}
EL^{q,p}_{\ell,R,S}H^1(X)=EL^{q,p}_{\ell,R,S}\bar{H}^1(X), 
\end{align*}
for all $p,q$ such that $1\leq p< d,$ $\frac{1}{p}-\frac{1}{q}=\frac{1}{\max\{n,d\}}$.
\end{pf}

\begin{cor}
\label{isop=>snp}
A bounded geometry Riemannian manifold which has isoperimetric dimension $\geq d>1$ is non-$(p,\ell,R,S)$-parabolic for all $1\leq p<d$ and all large enough $\ell$, $R$, $S\geq R$.
\end{cor}

\begin{pf}
This follows from Lemmata \ref{lpco=>snp} and \ref{isopdim=>redcoh}.
\end{pf}

\bigskip

It turns out that the isoperimetric dimensions of finitely generated groups are known.

\begin{prop}[Compare M. Troyanov, \cite{T}]
\label{isodimgps}
Let $G$ be a finitely generated group. Then the isoperimetric dimension of $G$ is
\begin{itemize}
  \item either equal to 1 if $G$ is virtually cyclic,
  \item or equal to its homogeneous dimension, an integer larger than 1 if $G$ is virtually nilpotent but not virtually cyclic.
  \item Otherwise, it is infinite.
\end{itemize}
If follows that a finitely generated group is $p$-parabolic if and only if it is virtually nilpotent of homogeneous dimension $\leq p$.
\end{prop}

\begin{pf}
According to T. Coulhon-L. Saloff Coste, \cite{CSC}, for finitely generated (or Lie) groups, volume growth provides an estimate on isoperimetric dimension. In particular, it implies that isoperimetric dimension is infinite unless volume growth is polynomial, in which case isoperimetric dimension is equal to the polynomial degree of volume growth. The only finitely generated groups of linear growth are virtually cyclic ones. That groups of polynomial growth are virtually nilpotent is M. Gromov's theorem of \cite{Gpol}. The isoperimetry of nilpotent groups was originally due to N. Varopoulos, \cite{VSC}.
\end{pf}

\subsection{Gr\"otzsch invariant}

Following \cite{grotzsch}, we use capacities to define a kind of large scale conformally invariant distance on a metric space.

\begin{defi}
Let $X$ be a metric space. Fix parameters $p,\ell,R,S$. For $x_1$, $x_2\in X$, let
\begin{align*}
\delta_{p,\ell,R,S,\phi}(x_1,x_2)&=&\inf\{\mathrm{cap}_{p,\ell,R,S}(\mathrm{im}(\gamma))\,;\\
&&\,\gamma\, \textrm{continuous arc in }X\textrm{ from }x_1\textrm{ to }x_2\}.
\end{align*}
\end{defi}

\begin{lem}
\label{deltainv}
Let $f:X\to X'$ be a large scale conformal map. Assume that $f$ is a bijection and that $f^{-1}:X'\to X$ is continuous. For all $R'>0$, there exists $R>0$ such that for all $\ell'\geq 1$, there exists $\ell\geq 1$ and $N'$ such that, for all $x_1$, $x_2\in X$ and all $S\geq R$,
\begin{align*}
\delta_{p,\ell,R,S}(x_1,x_2)\leq N'\,\delta_{p,\ell',R',\infty}(f(x_1),f(x_2)).
\end{align*}
\end{lem}

\begin{pf}
If $\gamma'$ is a continuous arc joining $f(x_1)$ to $f(x_2)$, $f^{-1}\circ\gamma'$ is a continuous arc joining $x_1$ to $x_2$, thus $\delta_{p,\ell,R,S}(x_1,x_2)\leq \mathrm{cap}_{p,\ell,R,\infty}(f^{-1}\circ\gamma')$.
Therefore $\delta_{p,\ell,R,\infty}(x_1,x_2)\leq N'\mathrm{cap}_{p,\ell,R,\infty}(\gamma')$, and taking an infimum,
\begin{align*}
\delta_{p,\ell,R,S}(x_1,x_2)\leq \delta_{p,\ell,R,\infty}(x_1,x_2)\leq N'\,\delta_{p,\ell',R',\infty}(f(x_1),f(x_2)).
\end{align*}
where the first inequality exploits the fact that adding constraints on packings decreases energies and capacities.
\end{pf}

\subsection{Upper bounds on capacities}

\begin{defi}
Say a metric space $X$ has \emph{controlled balls} if there exist $R>0$ and a measure $\mu$ and continuous functions $v>0$ and $V<\infty$ on $[R,+\infty)$  such that for every $x\in X$ and every $r\geq R$,
\begin{align*}
v(r)\leq\mu(B(x,r))\leq V(r).
\end{align*}
If such an estimate holds also for $r\in(0,R]$ and furthermore
\begin{align*}
\forall r\in(0,R], \quad v(r)\geq C\,r^Q,
\end{align*}
one says that $X$ is \emph{locally $Q$-Ahlfors regular}.
\end{defi}
In a Riemannian manifold or a simplicial complex with bounded geometry, balls are automatically controlled.

\begin{lem}
\label{upperPi}
Let $X$ be a geodesic metric space which has controlled balls. Let $p\geq 1$. Then $\delta_{p,\ell,R,\infty}$ is bounded above uniformly in terms of distance $d$. I.e. there exists a function $\Pi_{p,\ell,R}$ such that, if $d(x_1,x_2)\geq R$,
\begin{align*}
\delta_{p,\ell,R,\infty}(x_1,x_2)\leq \Pi_{p,\ell,R}(d(x_1,x_2)).
\end{align*}
If furthermore $X$ is locally $Q$-Ahlfors regular for some $Q\leq p$, then such an upper bound still holds with $R=0$, i.e. there exists a function $\Pi_{p,\ell}$ such that 
\begin{align*}
\delta_{p,\ell,0,\infty}(x_1,x_2)\leq \Pi_{p,\ell}(d(x_1,x_2)).
\end{align*}
\end{lem}

\begin{pf}
Fix $r\geq R$, set $T(r)=2+\frac{8r}{\ell-1}$ and
\begin{align*}
C(r)=\sup_{\rho\in[R,T]}\frac{(2\rho)^p}{v(\rho)}.
\end{align*}

For each $x\in X$, define a function $u_{x,r}$ as follows: $u_{x,r}=1$ on $B=B(x,r)$, vanishes outside $2B$ and is linear in the distance to $x$ in between. Let us estimate its $p$-energy. Let $\{B_j\}$ be a $(\ell,R,\infty)$-packing of $X$. If a ball $B_i$ intersects $2B$, and has radius $>4r/(\ell-1)$, then $\ell B_i$ contains $2B$. No other ball of the packing can intersect $2B$, hence an upper bound on $\sum_j\mathrm{diameter}(u_{x,r}(B_j))^p \leq 1$. Otherwise, all balls of the $\ell$-packing contributing to energy are contained in $TB$, for $T=2+(8r/(\ell-1)$. For each ball $B_j$ of radius $\rho_j$, $\rho_j\in [R,T]$ and
\begin{align*}
\mathrm{diameter}(u_{x,r}(B_j))^p
&\leq \mathrm{diameter}(B_j)^p\\
&\leq (2\rho_j)^p\\
&\leq C(r)\,v(\rho_j)\\
&\leq C(r)\,\mu(B_j).
\end{align*}
Summing up,
\begin{align*}
\sum_{j}\mathrm{diameter}(B_j)^p&\leq C(r)\,\sum_{j}\mu(B_j)\\
&\leq C(r)\,\mu(B(x,T))\\
&\leq C(r)V(T(r)).
\end{align*}
This gives an upper bound on $\mathrm{cap}_{p,\ell,R,\infty}(B(x,r))$ which depends on its radius $r$, on $p$, on $R$ and on $\ell$ only.

If $X$ is locally $Q$-Ahlfors regular and $Q\leq p$, then 
\begin{align*}
C'(r)=\sup_{\rho\in(0,T]}\frac{(2\rho)^p}{v(\rho)}<\infty,
\end{align*}
so the argument generalizes to arbitrary $(\ell,0,\infty)$-packings.

Balls in locally compact geodesic metric spaces contain geodesics which are continuous arcs. Thus the lower bound $\delta_{p,\ell,R,\infty}$ is bounded above by the capacity of a geodesic segment, which in turn is bounded above by the capacity of a ball, which is estimated in terms of its radius, on $p$ and on $\ell$ only.
\end{pf}

\subsection{Strong non-parabolicity}

Here, we are concerned with lower bounds on Gr\"otzsch' invariant $\delta$.

\begin{defi}
Let $X$ be a metric space. Say that $X$ is \emph{strongly non-$(p,\ell$, $R,S)$-parabolic} if $\delta_{p,\ell,R,S}(x_1,x_2)$ tends to infinity uniformly with $d(x_1,x_2)$. In other words, for every $R\leq S$ and $\ell>1$, there exists a function $\pi_{p,\ell,R,S}$ such that $\pi_{p,\ell,R,S}(T)$ tends to infinity when $T\to\infty$, and such that
\begin{align*}
\delta_{p,\ell,R,S}(x_1,x_2)\geq\pi_{p,\ell,R,S}(d(x_1,x_2)).
\end{align*}
\end{defi}

\begin{prop}
\label{lp=>snp}
Let $X$ be a metric space. Fix $\ell\geq 1$ and $S\geq R>0$. If $L^{q,p}_{\ell,R,S}H^1(X)=L^{q,p}_{\ell,R,S}\bar{H}^1(X)$ for some finite $q$, then $X$ is strongly non-$(p,\ell,R,S)$-parabolic.
\end{prop}

\begin{pf}
Let $x_1,x_2\in X$. Let $\gamma$ be a continuous arc joining $x_1$ to $x_2$. Fix $R>0$ and $\ell\geq 1$. Assume that $d(x_1,x_2)\geq 2R$. For each $j=0,\ldots,k:=\lfloor d(x_1,x_2)/2\ell R\rfloor$, pick a point $y_j$ on $\gamma$ such that $d(y_j,x_1)=2\ell Rj$. Let $B_j=B(y_j,R)$. By construction, $\{B_j\}$ is a $(\ell,R,R)$-packing of $X$. Let $u:X\to[0,1]$ be a function of bounded support such that $u=1$ on $\gamma$. Then $\sup_{B_j}u=1$, thus
\begin{align*}
\|u\|_{L^{q}_{\ell,R,R}}\geq k^{1/q}.
\end{align*}
By assumption, there exists $C$ such that, for every function $u$ of bounded support, 
\begin{align*}
\|u\|_{L^{q}_{\ell,R,R}}\leq C\,E^{p}_{\ell,R,R}(u)^{1/p}.
\end{align*}
This shows that
\begin{align*}
\mathrm{cap}_{p}(\gamma)\geq C\,k^{p/q}= C\,\lfloor\frac{d(x_1,x_2)}{2\ell R}\rfloor^{p/q},
\end{align*}
this is a lower bound on $\delta_{p,\ell,R,R}(x_1,x_2)$. This yields a lower bound on $\delta_{p,\ell,R,S}(x_1,x_2)$ for any $S$.\end{pf}

%\begin{prop}
%\label{lp=>snp}
%Let $X$ be a metric space which is uniformly perfect in the large. Fix $\ell\geq 1$ and $S\geq R>0$. If $L^{q,p}_{\ell,R,S}H^1(X)=L^{q,p}_{\ell,R,S}\bar{H}^1(X)$ for some finite $q$, then $X$ is strongly non-$(p,\ell,R,S)$-parabolic.
%\end{prop}

%\begin{pf}
%Let $x_1,x_2\in X$. Let $\gamma$ be a $\phi$-curve joining $x_1$ to $x_2$. For every $R>0$, the unit balls $B_j$ of $\llbracket 0,L_\gamma\rrbracket$ are mapped to balls $B'_j$ of radii $\geq R$. For every $\ell'\geq 1$, thery form a $(N',\ell',R,\infty)$-packing, for $N'=\phi_R(\ell')$.

\subsection{Consequences}

\begin{cor}
\label{snpp=>lcs=>qi}
Let $X$ and $X'$ be geodesic metric spaces which are strongly non-$(p,\ell,R,S)$-parabolic for some $p\geq 1$. Assume that both have controlled balls. Let $f:X\to X'$ be a homeomorphism such that both $f$ and $f^{-1}$ are large scale conformal maps. Then $f$ is a quasi-isometry.
\end{cor}

\begin{pf}
By strong non-$(p,\ell,R,S)$-parabolicity, $\delta$ invariants in $X$ are bounded below,
\begin{align*}
\delta_{p,\ell,R,S}(x_1,x_2)\geq \pi_{p,\ell,R,S}(d(x_1,x_2)).
\end{align*}
According to Lemma \ref{upperPi}, they are bounded above in $X'$,
\begin{align*}
\delta_{p,\ell,R,\infty}(x'_1,x'_2)\leq \Pi_{p,\ell,R}(d(x'_1,x'_2)).
\end{align*}
If $f:X\to X'$ is a large scale conformal homeomorphism, $N'\delta\circ f\geq\delta$ up to changes in parameters (Lemma \ref{deltainv}), 
\begin{align*}
N'\delta_{p,\ell',R',\infty}(f(x_1),f(x_2))\geq \delta_{p,\ell,R,S}(x_1,x_2).
\end{align*}
Combining these inequalities, we get for $f$,
\begin{align*}
N'\Pi_{p,\ell',R'}(d(f(x_1),f(x_2)))\geq \pi_{p,\ell,R,S}(d(x_1,x_2)),
\end{align*}
and for $f^{-1}$,
\begin{align*}
N'\Pi_{p,\ell',R'}(d(f^{-1}\circ f(x_1),f^{-1}\circ f(x_2)))\geq \pi_{p,\ell,R,S}(d(f(x_1),f(x_2))),
\end{align*}
hence
\begin{align*}
\pi_{p,\ell,R,S}(d(f(x_1),f(x_2)))\leq N'\Pi_{p,\ell',R'}(d(x_1,x_2)).
\end{align*}
These inequalities show that $f$ is a quasi-isometry.
\end{pf}

\begin{cor}
Let $M$, $M'$ be bounded geometry Riemannian manifolds quasi-isometric to non virtually cyclic groups. Then homeomorphisms $M\to M'$ which are large scale conformal in both directions must be quasi-isometries.
\end{cor}

\begin{rem}
This applies to Euclidean spaces of dimension $\geq 2$. Note that Examples \ref{power} are roughly conformal in both directions, but large scale conformal in only one direction.
\end{rem}

\begin{rem}
A natural question (Sylvain Maillot) is wether two spaces $X$ and $X'$ can have a large scale conformal map $X\to X'$ and a large scale conformal map $X'\to X$ without being quasi-isometric.
\end{rem}

\subsection{From coarse to uniformly conformal maps}

\begin{prop}
\label{snpp=>cc=>lcs}
Let $X$ be a metric space which is strongly non-$(p,\ell,R,S)$-parabolic for some $p\geq 1$ and all $\ell\geq 1$. Let $X'$ be a metric space which is locally $Q$-Ahlfors regular for some $Q\leq p$. 
Every coarsely conformal map $X\to X'$ is uniformly conformal.
Every roughly conformal map $X\to X'$ is large scale conformal.
\end{prop}

\begin{pf}
Fix $S\geq R>0$. Assume that $f:X\to X'$ is coarsely (resp. roughly) conformal. It suffices to show that for all $T'>0$, there exists $T_0>0$ such that $f$ maps no $T$-ball $B$, $T\geq T_0$, into a $T'$-ball $B'$. Given $\ell'\geq 1$, there exist $\ell\geq 1$ and $N'$ such that if $f(B)\subset B'$, $\mathrm{cap}_{p,\ell,R,S}(B)\leq N'\mathrm{cap}_{p,\ell',0,\infty}(B')$. This upper bound fails if $T$ is sufficiently large, $T\geq T_0$. This shows that $T$-balls are never mapped into $T'$-balls. So if $f$ is coarsely conformal, it is in fact $(T_0,S_0,T',\infty)$-coarsely conformal, for large enough $T_0$ and for every $S_0\geq T_0$, hence it is uniformly conformal. If $f$ is roughly conformal, it is in fact $(T_0,\infty,T',\infty)$-coarsely conformal, for large enough $T_0$ thus $f$ is large scale conformal.
\end{pf}

\begin{rem}
The assumptions of Proposition \ref{snpp=>cc=>lcs} are satisfied for $X=\R^n$ provided $p<n$ and for $X'=\R^{n'}$ for $p\geq n'$. So Proposition \ref{snpp=>cc=>lcs} applies if $n'<n$, i.e. exactly when there are no coarse conformal maps $X\to X'$. In fact, the conclusion fails if $n=n'$, as Examples \ref{power} show.
\end{rem}

Proposition \ref{snpp=>cc=>lcs} allows to modify the assumptions in the corollaries of subsections \ref{tonil} and \ref{tohyp}. For instance,

\begin{cor}\label{arbtonil}
Let $G$ be a finitely generated group. Let $G'$ be a nilpotent Lie group equipped with a left-invariant Riemannian metric. If there exists a coarse conformal map $G\to G'$, then $G$ is virtually nilpotent and $d(G)\leq d(G')$.
\end{cor}

\begin{pf}
If $G$ is virtually cyclic, then it is virtually nilpotent and $d(G)=1\leq d(G')$. Otherwise, $G$ has isoperimetric dimension $Q>1$, thus it is strongly non-$(p,\ell,R,S)$-parabolic for all $1<p<Q$. $G'$ is locally $n'$-Ahlfors regular for $n'=\mathrm{dimension}(G')$. Assume that $Q>Q'$. Since $n'\leq Q'$, one can pick $p$ such that $\max\{1,n'\}< p<Q$. Proposition \ref{snpp=>cc=>lcs} asserts that a coarse conformal map $G\to G'$ is automatically uniformly conformal. Corollary \ref{niltonil} shows that such a map cannot exist.\end{pf}

\begin{cor}
\label{hyptohypp}
Let $G$, $G'$ be non-elementary hyperbolic groups. Let $M'$ be a Riemannian manifold of bounded geometry, which is quasi-isometric to $G'$. If there exists a coarse conformal map $G\to M'$, then 
$$\mathrm{CohDim}(G)\leq \mathrm{ConfDim}(G').$$
\end{cor}

\begin{pf}
Non-elementary hyperbolic groups have infinite isoperimetric dimensions. Thus $G$ is strongly non-$(p,\ell,R,S)$-parabolic for all $p$. By assumption, $M'$ is locally $n'$-Ahlfors regular for $n'=\mathrm{dimension}(M')$. Choose some $p\geq n'$. Proposition \ref{snpp=>cc=>lcs} asserts that a coarse conformal map $G\to M'$ is automatically uniformly conformal. Composing with a quasi-isometry, we get a uniformly conformal map $G\to G'$, so Corollary \ref{hyptohyp} applies.
\end{pf}

Corollary \ref{cdd} is a combination of Corollaries \ref{arbtonil} and \ref{hyptohypp}.

\subsection{Large scale conformality in one dimension}

We have been unable to extend Theorem \ref{isomorphism} to the virtually cyclic case. Here is a partial result.

\begin{lem}
Let $f$ and $g$ be continuous maps $\R\to\R$ such that 
\begin{itemize}
  \item $f$ and $g$ are large scale conformal;
  \item $g\circ f$ and $f\circ g$ are coarse embeddings.
\end{itemize}
Then $f$ is a quasi-isometry.
\end{lem}

\begin{pf}
Fix $R'\leq R$, $\ell'$, $\ell$, $N'$ as given by the definition of large scale conformality. Let $\tilde{R}$ be given by the definition of coarse embeddings: $g\circ f$ maps $R$-balls to $\tilde{R}$-balls. To save notation, assume that the same constants serve for $g$. Since we are on the real line, balls of radius $R$ are intervals of length $2R$. In the correspondence between balls $B\mapsto B'$, one can assume that $B'$ is a minimal interval containing $f(B)$ and of length $\geq 2R'$, i.e. $B'=f(B)$ itself if length$(f(B))\geq 2R'$.

The balls $B_j=B(2\ell Rj,R)$ are mapped into balls $B'_j$ forming an $(N',\ell',R',$ $\infty)$-packing. Assume that $f(B_0)$ has length $2R'_0\geq 2R'$, in order that $B'_0=f(B_0)$. Let $\{B''_j\}$ be an $(\ell,R,\infty)$-packing of $B'_0$. The number of balls in this packing is at least $\frac{R'_0}{2\ell R}$. In turn, $g$ maps $B''_j$ into $B'''_j$, which form an $(N',\ell',R',\infty)$-packing, which is the union of $N'$ $(\ell',R',\infty)$-packings. One of them has at least $\frac{R'_0}{2\ell RN'}$ elements. Every ball $B'''_j$ contains 
$$g(B''_j)\subset g(B'_0)=g\circ f(B_0)\subset \tilde{B}:=B(g\circ f(0),\tilde{R}),$$
so $B'''_j$ intersects $\tilde{B}$. At most two of these balls contains boundary points, so all others are contained in $\tilde{B}$. At least one of these balls has radius $\rho$ no larger than
\begin{align*}
\frac{\tilde{R}}{\frac{R'_0}{2\ell RN'}-2}=\frac{2\ell RN'\tilde{R}}{R'_0-4\ell RN'}.
\end{align*}
Since $\rho\geq R'$, we obtain an upper bound on $R'_0$. This shows that all $R$-balls are mapped to $R'_0$-balls, i.e. $f$ is a coarse embedding.
\end{pf}

\bigskip

\noindent
Pierre Pansu 
\par\noindent Laboratoire de Math\'ematiques d'Orsay,
\par\noindent Univ. Paris-Sud, CNRS, Universit\'e
Paris-Saclay
\par\noindent 91405 Orsay, France.
\par\noindent
e-mail: pierre.pansu@math.u-psud.fr

\bigskip

Key words: conformal; quasiconformal; metric space; Dirichlet energy; harmonic; cohomology; capacity.

\bigskip

MSC: 53A30; 53C23; 52C26; 33E05; 31C12; 31A15; 30F45; 37F35.

\begin{thebibliography}{}

\bibitem[A]{A} Patrice Assouad,
{\em Plongements lipschitziens dans $\R^{n}$}. 
Bull. Soc. Math. France \textbf{111} (1983), 429--448.

\bibitem[B]{B} Marc Bourdon,
{\em Une caract\'erisation alg\'ebrique des hom\'eomorphismes quasi-M\"obius.}
Ann. Acad. Sci. Fenn. Math.  \textbf{32}  (2007), 235--250.

\bibitem[BC]{BC} Itai Benjamini; Nicolas Curien,
\emph{On limits of graphs sphere packed in Euclidean space and applications}.
European J. Combin. \textbf{32} (2011), no. 7, 975--984.

\bibitem[BK]{BK} Mario Bonk; Bruce Kleiner,
\emph{Conformal dimension and Gromov hyperbolic groups with $2$-sphere boundary.}
Geom. Topol. \textbf{9} (2005), 219--246.

\bibitem[BnK]{BnK} Marc Bourdon; Bruce Kleiner,
\emph{Some applications of $l_p$-cohomology to boundaries of Gromov hyperbolic spaces}.
Groups Geom. Dyn. \textbf{9} (2015), no. 2, 435--478.
arXiv:1203.1233

%\bibitem[Bn]{Bn} Mario Bonk,
%{\em Quasiconformal geometry of fractals.} 
%International Congress of Mathematicians. Vol. II, 1349--1373, Eur. Math. Soc., Z\"urich, 2006.

\bibitem[Bo]{Bo} Marc Bourdon,
{\em Sur les immeubles fuchsiens et leur type de quasi-isom\'etrie.}
Erg. Th. Dynam. Syst. \textbf{20} (2000), 343--364. 

%\bibitem[BHT]{BHT} Zolt\'an Balogh, Ilkka Holopainen and Jeremy Tyson,
%{\sl Singular solutions, homogeneous norms, and quasiconformal mappings in Carnot groups. }
%Math. Ann. \textbf{324} (2002), 159--186. 

\bibitem[BM]{Berard-Meyer} Pierre B\'erard; Daniel Meyer,
\emph{In\'egalit\'es isop\'erim\'etriques et applications.}
Ann. Sci. \'Ecole Norm. Sup. (4) \textbf{15} (1982), 513--541.

\bibitem[BP]{BP} Marc Bourdon; Herv\'e Pajot, 
{\sl Cohomologie $l_p$ et espaces de Besov.} 
J. Reine Angew. Math. \textbf{558} (2003), 85--108.

\bibitem[BS]{BS} Itai Benjamini; Oded Schramm,
\emph{Lack of sphere packing of graphs via non-linear potential theory}.  
J. Topol. Anal. 5 (2013), no. 1, 1--11.

%\bibitem[BST]{BST} Itai Benjamini, Oded Schramm and \'Adam Tim\'ar,
%\emph{On the separation profile of infinite graphs}.  
%Groups Geom. Dyn. \textbf{6} (2012), no. 4, 639--658.

\bibitem[BoS]{BoS} Mario Bonk; Oded Schramm,
\emph{Embeddings of Gromov hyperbolic spaces.}
Geom. Funct. Anal. \textbf{10} (2000), 266--306. 

\bibitem[Ca]{carrasco} Matias Carrasco-Piaggio,
\emph{On the conformal gauge of a compact metric space.}
Ann. Sci. \'Ecole Norm. Sup. \textbf{46}, no 3 (2013), 495--548.

\bibitem[C]{Coo} Michel Coornaert, 
\emph{Mesures de Patterson-Sullivan sur le bord d'un espace hyperbolique au sens de Gromov}. 
Pacific J. Math. \textbf{159} (1993), no. 2, 241--270.

\bibitem[CSC]{CSC} Thierry Coulhon; Laurent Saloff-Coste,
{\em Isop\'erim\'etrie pour les groupes et les vari\'et\'es.} 
Rev. Mat. Iberoamericana \textbf{9} (1993), 293--314.

\bibitem[CT]{CT} Yves Cornulier; Romain Tessera,
{\em Contracting automorphisms and $L^p$-cohomology in degree one}.
Arkiv Mat. \textbf{49:2} (2011), 295--324. \texttt{arXiv:0908.2603}.

\bibitem[D]{ducret} Stephen Ducret,
\emph{$L^{q,p}$-Cohomology of Riemannian Manifolds
and Simplicial Complexes of Bounded Geometry.}
Th\`ese n$^o$ 4544, Ec. Pol. Fed. Lausanne (2009).

\bibitem[E]{E} G\'abor Elek, 
{\sl Coarse cohomology and $l_p$-cohomology.} 
$K$-Theory \textbf{13} (1998), 1--22.

%\bibitem[F]{Federer} Herbert Federer,
%{\sl Geometric measure theory}. 
%Die Grundlehren der mathematischen Wissenschaften, Band \textbf{153} Springer-Verlag New York Inc., New York (1969).

\bibitem[G1]{Gpol} Mikhael Gromov,
\emph{Groups of polynomial growth and expanding maps.}
Inst. Hautes \'Etudes Sci. Publ. Math. \textbf{53} (1981), 53--73.

\bibitem[G2]{G} Mikhael Gromov,
\emph{Asymptotic invariants of infinite groups}. 
Geometric group theory, Vol. 2 (Sussex, 1991), 1--295, London Math. Soc. Lecture Note Ser., \textbf{182}, Cambridge Univ. Press, Cambridge, 1993.

\bibitem[Ge]{genton} Luc Genton,
\emph{Scaled Alexander-Spanier Cohomology and $L^{q,p}$ Cohomology for Metric Spaces.}
Th\`ese n$^o$ 6330, Ec. Pol. Fed. Lausanne (2014).

\bibitem[GH]{GH} Etienne Ghys; Pierre de la Harpe, 
\emph{Sur les groupes hyperboliques d'apr\`es Mikhael Gromov} (Bern, 1988), 1--25, Progr. Math., \textbf{83}, Birkh\"auser Boston, Boston, MA, 1990.

\bibitem[Gr]{grotzsch} Herbert Gr\"otzsch,
\emph{\"Uber die Verzerrung bei schlichten nichtkonformen Abbildungen und \"uber eine damit Zusammenh\"angende Erweiterung des Picardschen Satzes.}
Ber. Verh. S\"achs. Akad. Wiss. Leipzig \textbf{80} (1928), 503-507.

\bibitem[GT]{GT} Vladimir Gol'dshtein; Marc Troyanov,
\emph{Axiomatic Sobolev spaces on metric spaces.} 
Function spaces, interpolation theory and related topics (Lund, 2000), 333--343, de Gruyter, Berlin, 2002.

\bibitem[H]{H} Ernst Heintze,
\emph{On homogeneous manifolds of negative curvature}. 
Math. Ann. \textbf{211} (1974), 23--34.

%\bibitem[H]{H} Ursula Hamenst\"adt, 
%\emph{A geometric characterization of negatively curved locally symmetric spaces.}
%J. Differential Geom. \textbf{34} (1991), 193--221.

%\bibitem[HaP]{HaP} Peter Ha\"{\i}ssinsky and Kevin Pilgrim,
%{\em Coarse expanding conformal dynamics.} 
%Ast\'erisque \textbf{325} (2009), 1--139.

\bibitem[He]{He} Juha Heinonen, 
\emph{Lectures on analysis on metric spaces.} 
Universitext. Springer-Verlag, New York, 2001.

\bibitem[HK]{HK} Juha Heinonen; Pekka Koskela,
\emph{Definitions of quasiconformality.} 
Invent. Math. \textbf{120} (1995), 61--79. 

\bibitem[HMT]{HMT} David Hume; John M. Mackay; Romain Tessera,
\emph{Poincar\'e profiles of groups and spaces.}
arXiv:1707.02151 

%\bibitem[HP]{HP} Sa'ar Hersonsky and Fr\'ed\'eric Paulin, 
%{\em On the rigidity of discrete isometry groups of negatively curved spaces.} 
%Comment. Math. Helv. \textbf{72} (1997), 349--388. 

%\bibitem[LS]{LS} Jouni Luukkainen and Eero Saksman,
%\emph{Every complete doubling metric space carries a doubling measure}.
%Proc. Amer. Math. Soc. \textbf{126} (1998), 531--534. 

\bibitem[MAG]{metric2011} Collective blog,

https://metric2011.wordpress.com/2011/01/24/notes-of-james-lees-lecture-nr-1/

\bibitem[McK]{McK} John Mackay,
\emph{Conformal dimension via subcomplexes for small cancellation and random groups.}
Math. Annalen. \textbf{364:3--4} (2016), 937--982. arxiv:1409.0802.

\bibitem[MN]{Munoz-Flores Nardulli} Abraham Munoz-Flores; Stefano Nardulli,
\emph{Continuity and differentiability properties of the isoperimetric profile in complete noncompact Riemannian manifolds with bounded geometry.}
http://cvgmt.sns.it/paper/2408/

\bibitem[NN]{NN} Assaf Naor; Ofer Neiman,
\emph{Assouad's theorem with dimension independent of the snowflaking}. 
Rev. Mat. Iberoam. \textbf{28} (2012), no. 4, 1123--1142. 

%\bibitem[P]{P} Samuel S. Patterson, 
%{\em The limit set of a Fuchsian group.} 
%Acta Math. \textbf{136} (1976), 241--273.

\bibitem[P1]{growth} Pierre Pansu, 
{\em Croissances des boules et des g\'eod\'esiques ferm\'ees dans les nilvari\'et\'es.}
Ergod. Th. Dynam. Syst. \textbf{3}, (1983), 415--445.

\bibitem[P2]{Pa} Pierre Pansu, 
{\em Dimension conforme et sph\`ere \`a l'infini des vari\'et\'es \`a courbure n\'egative.}
Ann. Acad. Sci. Fenn. Ser. A I Math. \textbf{14} (1989), 177--212. 

\bibitem[P3]{Pa2} Pierre Pansu, 
{\em Cohomologie $L^p$ des vari\'et\'es \`a courbure n\'egative, cas du degr\'e 1.}
Rend. Semin. Mat. Torino, Fasc. Spec.,  (1989), 95--120.

\bibitem[P4]{Pa1} Pierre Pansu, 
{\em Cohomologie $L^p$ en degr\'e 1 des espaces homog\`enes.}
Potential Anal. \textbf{27} (2007), 151--165.

\bibitem[P5]{Pa4} Pierre Pansu, 
{\em Cohomologie $L^p$ et pincement.}
Comment. Math. Helv. 83 (2008), 327--357.

\bibitem[P6]{BFP} Pierre Pansu,
\emph{Cup-products in $L^{q,p}$-cohomology: discretization and quasi-isometry invariance.}
arXiv:1702.04984.

%\bibitem[Rb]{Rb} Thomas Roblin,
%{\em Ergodicit\'e et \'equidistribution en courbure n\'egative.}
%M\'em. Soc. Math. Fr. \textbf{95} (2003). 

%\bibitem[Re]{Re} John Roe,
%{\em Coarse cohomology and index theory on complete Riemannian manifolds.} 
%Mem. Amer. Math. Soc. \textbf{104} (1993).

\bibitem[R]{R} Halsey Royden, 
{\em On the ideal boundary of a Riemann surface.}  
Contributions to the theory of Riemann surfaces,  pp. 107--109. Annals of Mathematics Studies \textbf{30}. Princeton University Press, Princeton, N. J., 1953.

\bibitem[RS]{RS} Burt Rodin; Dennis Sullivan,
{\sl The convergence of circle packings to the Riemann mapping.}
J. Differential Geom. \textbf{26} (1987), 349--360. 

\bibitem[S]{S} Dennis Sullivan,
{\sl The density at infinity of a discrete group of hyperbolic motions.}
Inst. Hautes Etudes Sci. Publ. Math. \textbf{50} (1979), 171--202. 

\bibitem[Sh]{Shchur} Vladimir Shchur,
{\sl On the quantitative quasiisometry
problem: Transport of Poincar\'e
inequalities and different types of quasiisometric
distortion growth.} 
J. Funct. Anal. \textbf{269} (2015), 3147--3194.

\bibitem[ST]{ST} Daniel A. Spielman; Shang-Hua Teng, 
{\em Disk packings and planar separators.}
12th Symp. Computational Geometry, (1996), 349--359.

\bibitem[T]{T} Marc Troyanov,
{\sl Parabolicity of manifolds.}
Siberian Adv. Math. \textbf{9} (1999), 125--150. 

\bibitem[Ty]{Ty} Jeremy Tyson,
{\sl Metric and geometric quasiconformality in Ahlfors regular Loewner spaces.} 
Conform. Geom. Dyn. \textbf{5} (2001), 21--73 (electronic).

\bibitem[V1]{vaisala} Jussi V\"ais\"al\"a,  
\emph{Lectures on $n$-dimensional quasiconformal mappings.}
Lecture Notes in Mathematics, \textbf{229}. Springer-Verlag, Berlin-New York, (1971).

\bibitem[V2]{V} Jussi V\"ais\"al\"a,  
{\em Quasi-M\"obius maps.} 
J. Analyse Math. \textbf{44} (1984/85), 218--234.

%\bibitem[VK]{VK} Alexander L. Vol'berg and Sergei V. Konyagin, 
%{\em On measures with the doubling condition.} 
%Izv. Akad. Nauk SSSR Ser. Mat. \textbf{51} (1987), no. 3, 666--675; translation in Math. USSR-Izv. \textbf{30} (1988), no. 3, 629--638.

\bibitem[VSC]{VSC} Nicholas Varopoulos; Laurent Saloff-Coste; Thierry Coulhon,
\textrm{Analysis and geometry on groups}. 
Cambridge Tracts in Mathematics, 100. Cambridge University Press, Cambridge, 1992.

\end{thebibliography}
\end{document}